\documentclass[a4paper]{article} 

\usepackage{latexsym}
\usepackage{amsmath}
\usepackage{amsfonts}
\usepackage{amssymb}
\usepackage{graphicx}
\usepackage{color}
\usepackage{epsfig}

\def\leq{\leqslant}
\def\geq{\geqslant}

\def\N{\mathbb{N}}
\def\R{\mathbb{R}}
\def\Z{\mathbb{Z}}
\def\Q{\mathbb{Q}}

\newtheorem{Pro}{Proposition}[section]
\newtheorem{Thm}[Pro]{Theorem}
\newtheorem{Lem}[Pro]{Lemma}

\begin{document}

\title{Discrete time piecewise affine models of genetic regulatory networks} 
\author{R.\ Coutinho$^1$, B.\ Fernandez$^2$, R.\ Lima$^2$ and A.\
Meyroneinc$^2$}  
\maketitle
\begin{center}
$^1$ Departamento de Matem\'atica\\
Instituto Superior T\'ecnico\\
Av. Rovisco Pais 1096\\
Lisboa Codex Portugal\\
{\tt Ricardo.Coutinho@math.ist.utl.pt}
\bigskip

$^2$
Centre de Physique Th\'eorique\\
CNRS - Universit\'es de Marseille I et II, et de Toulon\\
Luminy Case 907\\ 
13288 Marseille CEDEX 09 France\\
{\tt fernandez@cpt.univ-mrs.fr, lima@cpt.univ-mrs.fr, meyroneinc@cpt.univ-mrs.fr}
\bigskip

{\bf Abstract}
\end{center} 
We introduce simple models of genetic regulatory networks and we proceed to the mathematical analysis of their dynamics. The models are discrete time dynamical systems generated by piecewise affine contracting mappings whose variables represent gene expression levels. When compared to other models of regulatory networks, these models have an additional parameter which is identified as quantifying interaction delays. In spite of their simplicity, their dynamics presents a rich variety of behaviours. This phenomenology is not limited to piecewise affine model but extends to smooth nonlinear discrete time models of regulatory networks. 

\noindent
In a first step, our analysis concerns general properties of networks on arbitrary graphs (characterisation of the attractor, symbolic dynamics, Lyapunov stability, structural stability, symmetries, etc). In a second step, focus is made on simple circuits for which the attractor and its changes with parameters are described. In the negative circuit of 2 genes, a thorough study is presented which concern  stable (quasi-)periodic oscillations governed by rotations on the unit circle -- with a rotation number depending continuously and monotonically on threshold parameters. These regular oscillations exist in negative circuits with arbitrary number of genes where they are most likely to be observed in genetic systems with non-negligible delay effects.

\bigskip
\section{Introduction}
With genome sequencing becoming a widespread procedure, the structural
information contained in the genome is now largely accessible. Not
much can be said about the functional information contained in gene
expression regulatory mechanisms which is still largely
unravelled. Under a simplifying point of view, regulatory mechanisms
are described in terms of networks of basic process in
competition. Insights on the role of regulatory networks during
the development of an organism and changes with environmental parameters can be
gained from the analysis of crude dynamical models \cite{TCN03}. The
models are usually supported by directed graphs where the nodes
represent genes (or their products) and where the arrows represent
interactions between genes.  

\noindent
Depending on context, various formalisms have been used to model regulatory networks, see \cite{D02} for a recent review. Discrete
variable models (boolean networks) have been employed to obtain
essential features, such as influence of the sign of a circuit on
multi-stationarity or homeostasis \cite{T73,TA90} or decomposition into
regulatory modules \cite{TR99}. In addition to the analysis of
behaviours for arbitrary graphs, models of specific regulatory
mechanisms have been derived and thoroughly investigated in this
formalism, e.g.\ flower morphogenesis in {\sl Arabidopsis thaliana}
\cite{MTA99} and dorso-ventral patterning in {\sl Drosophila
  melanogaster} \cite{SHT97}.

\noindent
In a more traditional framework, (systems of coupled) nonlinear
ordinary differential equations with interactions represented by
sigmoid functions and early piecewise affine analogues have been
considered \cite{GP78}.  Piecewise affine models admit analytical
investigation without affecting most features of the dynamics. These
models have been studied by using tools inspired from the theory of
dynamical systems. In \cite{E00}, a method has been developed to
determine existence and stability of periodic trajectories with
prescribed qualitative behaviour. In \cite{DGHPSG04} 
and in \cite{ESAG01}, symbolic dynamics has been employed in a
computational framework in order to obtain results on qualitative
behaviours and their changes with parameters, namely bifurcations. 
Naturally coupled differential equations have not only been analysed
in their own but they have also been applied to represent specific
mechanisms, see \cite{TCN03} for a review. 
\bigskip

In spite of being governed by the same rules (described below), boolean networks and
coupled ordinary differential equations often present distinct dynamical
behaviours \cite{TA90}. Boolean networks reveal periodic orbits where
differential equations only possess stationary points.

\noindent
These distinct behaviours however can be recovered in a unique and simple model
by adjusting a parameter, say $a$. The model, an original discrete
time dynamical system with continuous variables, obeys the same rules
as in previous models. Gene product concentrations (expression
levels) evolve according to combined interactions from other genes
in the network. The interactions are given by step functions which
express that a gene acts on another gene, or becomes inactive, only
when its product concentration exceeds a threshold. 

\noindent
Based on comparisons with (systems of coupled) delay differential equations, the parameter $a$ appears to be related with a time delay. When this parameter is such that the delay is maximum, the system reduces to a boolean network. And in the limit of a system without delay, the dynamics is as for a differential equation. It is therefore of particular interest to describe the phase portraits not only depending on interaction parameters but depending on the delay parameter $a$.

\noindent
That a discrete time dynamical system provides a simple substitute to a delay differential equation is a well-known fact which has already been employed in mathematical biology \cite{E96,KL93}. Note that this does not contradicts the statements on equivalence between discrete time dynamical systems and ordinary differential equations (Poincar\'e section). Indeed, such statements match the dynamics of discrete time dynamical systems in an $N$-dimensional phase space with the dynamics of an ordinary differential equations in an $N+1$-dimensional phase space.
\bigskip

The paper is organised as follows. With the definition of the model
provided (section \ref{R-MODEL}), we proceed to a comprehensive
analysis of its dynamics with emphasis on changes with parameters. The
analysis begins with the study 
of general properties of networks on arbitrary graphs and more
specifically on circuits (section \ref{R-GENER}). In 
particular, attractors are described in terms of symbolic dynamics by
means of an admissibility condition on symbolic sequences. Lyapunov
stability, structural stability and symmetries of orbits in the
attractor are presented. A special section discusses relationships
with boolean networks and differential equations. 

\noindent
In a second step, we focus on simplest feedback circuits with 1 and 2
genes (section \ref{R-SIMPL}). An analysis of
the dynamics is presented which is complete in phase space and in
parameter space for most of these circuits. The results rely on 
previous results on piecewise affine and contracting rotations.

\noindent
The dynamics of negative circuit with two genes requires special
attention and a fully original analysis. The most important orbits (the most likely to be seen in numerical experiments with systematic prospection in phase space) are the so-called regular orbits. These are stable (quasi-)periodic oscillations governed by rotations on the unit circle composed of arcs associated with atoms in phase space. The associated characteristics (rotation number and arc lengths) depend on systems parameters in agreement with the parameters dependent sojourn times in atoms of such oscillations.
\bigskip

For the sake of clarity, the paper is decomposed into two parts. In
part A (sections \ref{R-GENER} and \ref{R-SIMPL}), we only present
results. Some results can be obtained easily and their proof are left to the reader. Most original results however require elaborated mathematical proofs and calculations. These are postponed to part B (sections \ref{M-ANAL} to \ref{M-STABI}). 

\noindent
The final section \ref{M-LONG} contains concluding remarks and open problems. In particular, it indicates how the regular orbit analysis developed in the negative circuit of two genes extends naturally to circuits with arbitrary number of genes.
 
\section{The model}\label{R-MODEL}
Basic models of genetic regulatory mechanisms are networks of
interacting genes where each gene is submitted both to a
self-degradation and to interactions from other genes. The
self-degradation has constant rate and the interactions are
linear combinations of sigmoid functions. In discrete time, a simple
piecewise affine model can be defined by the following relation 
\begin{equation}
x_i^{t+1}=ax_i^t+(1-a)\sum_{j\in I(i)} K_{ij}
H(s_{ij}(x_j^t-T_{ij})),\quad i=1,N
\label{MODEL}
\end{equation}
i.e.\ $x^{t+1}=F(x^t)$ where $x^t=\{x_i^t\}_{i=1,N}$ is the local
variable vector at time $t\in\Z$ and
$F_i(x)=ax_i+(1-a){\displaystyle\sum_{j\in I(i)}} K_{ij}
H(s_{ij}(x_j-T_{ij}))$ is the $i$th component of the mapping $F$
defined from the phase space $\R^N$ into itself. 

\noindent
In expression (\ref{MODEL}), the subscript $i$ labels
a gene ($N$ denotes the number of genes involved in the network). The
graph supporting the network (the arrows between genes) is implicitly
given by the sets $I(i)\subset\{1,\cdots,N\}$. For each $i$ the set
$I(i)$ consists of the set of genes which have an action on $i$. In
particular, a self-interaction (loop) occurs when some set $I(i)$
contains $i$. Examples of networks are given Figure \ref{FIGCIR} and
Figure \ref{FIGFLIP}.  

\noindent
The property of the action from $j\in I(i)$ to $i$ is specified by a
sign, the number $s_{ij}$. An activation is associated with a positive
sign, i.e.\ $s_{ij}=+1$, and an inhibition with a
negative sign $s_{ij}=-1$. 

\noindent
The degradation rate $a\in [0,1)$ is supposed to be identical for all genes. This assumption only serves to simplify calculations and the resulting expression of existence domain of given orbits (bifurcation values). However, the whole analysis in the paper does not depend on this assumption and extends immediately to the case where the degradation rate depends on $i$.

\noindent
The symbol $H$ denotes the Heaviside function 
\[
H(x)=\left\{\begin{array}{l}
0\text{ if }x<0\\
1\text{ if }x\geq 0
\end{array}\right.
\]
In order to comply with the assumption of cumulative interactions, the
interaction intensities $K_{ij}$ are supposed to be positive. They are normalised
as follows 
\[
\sum_{j\in I(i)} K_{ij}=1,\quad i=1,N
\]
This normalisation is arbitrary and has no consequences on the dynamics (see section \ref{M-NORMAL} for a discussion on normalisation). These
properties imply that the 
(hyper)cube $[0,1]^N$ is invariant and absorbs the orbit of every
initial condition in $\R^N$. In another words, every local variable
$x_i^t$ asymptotically (when $t\to\infty$) belongs to the interval
$[0,1]$ (see section \ref{M-ATTRA}).  Ignoring the behaviour outside
$[0,1]^N$,  we may only consider the dynamics of initial conditions in
this set. When normalised to $[0,1]$, the variable $x_i^t$ should be
interpreted as a ratio of gene product concentration produced by the regulatory process, rather than as a chemical concentration. Lastly, the parameters $T_{ij}$ belong to the (open)
interval $(0,1)$ and represent interaction thresholds (see section
\ref{M-NORMAL} for a discussion on interaction thresholds domains).  

\vskip 1truecm

\centerline{\sc \large Part A. Results}

\section{General properties of the dynamics}\label{R-GENER}
In this section we present some dynamical properties of networks on
arbitrary graphs, and more restrictively, on arbitrary circuits. Most
properties are preliminaries results which allow us to simplify the
analysis of circuit dynamics to follow.  The mathematical analysis of
the results in this section is given in section \ref{M-ANAL}. 

\subsection{Symbolic dynamics of genetic regulatory networks}\label{R-ATTRA}
Following a widespread technique in the theory of dynamical systems, the
qualitative features (the structure) of a dynamical system can be 
described by using symbolic dynamics \cite{R99}. This consists in
associating sequences of symbols with orbits. To that goal, a coding
needs to be introduced which associates a label with each domain in
phase space. In our case of piecewise affine mapping,
the atoms are domains -- bounded by discontinuity lines -- where $F$ is affine. They are naturally labelled by the (elementary) symbols $\theta_{ij}=H(s_{ij}(x_j-T_{ij}))\in\{0,1\}$
involved in interactions. In
particular, the symbols \footnote{We also use the term symbol for the concatenation $(\theta_{ij})_{i=1,N,j\in I(i)}$ of elementary symbols.} 
depend on the number of genes, on the interaction graph and on the interaction signs. 

\noindent
By evaluating the atom(s) the image by $F$ of a given atom intersects,
a {\bf symbolic graph} is obtained which indicates the possible (one-step)
transitions between symbols. As a consequence every point in phase space
generates via its orbit, a symbolic sequence, its {\bf
  code}, which corresponds to an infinite path in the symbolic graph. 

\noindent
It may happen that a code is associated with two distinct
points (absence of injectivity). For instance, the code associated with an initial condition in the
immediate basin of attraction of a periodic point is the same as the
periodic point code. Simple examples can be found for
the self-activator, see section \ref{R-SELFACT}.

\noindent
It may also happen that an infinite path in the symbolic graph does
not correspond to any point in phase space. In the present framework of piecewise contracting mapping, this happens when there exists an atom whose image intersects several atoms (absence of Markov property).
The simplest such example is the self-inhibitor, see section \ref{R-CONTRAC}.
\bigskip

The injectivity of the coding map associated with $F$ can be shown to
hold in the {\bf attractor}. (The attractor is the set of points which
attracts all orbits in 
phase space, see section \ref{M-ATTRA} for a 
definition.) Considering the attractor amounts to focusing on
asymptotic dynamics. In applications, this is particularly relevant when transients are short.

\noindent
That distinct points in the attractor have distinct codes is a
consequence of the following property. Points in the attractor are completely determined by their code (and the parameters). The expression of their
coordinates is a 
uniformly converging series - relation (\ref{GLOBORB}) in section
\ref{M-ATTRA}. \footnote{Points $x^0$ in the attractor turn out to
  have pre-images $x^{-t}$ for any $t\in\N$ (see Proposition
  \ref{AEQUG}). Therefore, the symbols 
  $\theta^t_{ij}$ associated with points in the attractor are
  defined for all $t\in\Z$ and not only for all $t\in\N$, see relation
  (\ref{GLOBORB}).}

\noindent
In addition the relation (\ref{GLOBORB}) also provides a criterion for a 
symbolic sequence to code for a point in the attractor. (If it does, the symbolic sequence is said to be admissible.) The criterion, the {\bf
  admissibility condition} -- relation (\ref{ADMICOND}) --  simply
imposes that, for each $t$, the formal point 
$x^t$ computed by using relation (\ref{GLOBORB}) belongs to the atom
labelled by the symbols $\theta_{ij}^t$ and is therefore a genuine
orbit point.
\bigskip

In practise, the analysis of a genetic regulatory network consists in
analysing the corresponding admissibility condition (based on the
transition graph) in order to determine which symbolic sequences are admissible,
possibly depending on parameters. 
Note that according to relation (\ref{ADMICOND}), when intersected with any hyperplane $a=\text{constant}$, the admissibility domain of any given symbolic sequence reduces to a product of
threshold parameter intervals.

\subsection{Lyapunov stability and robustness with respect to changes
  in parameters}\label{R-RELEV} 
The assumption $a<1$, which reflects self-degradation of genes implies that $F$ is a piecewise contraction. Orbits of piecewise contractions are robust with
respect to changes in initial conditions and to changes in 
parameters. Such robustness have been identified as generic features
in various genetic regulatory networks, see \cite{MDMO02} and
references therein. 

The robustness properties concern orbits, in the attractor of $F$, not
intersecting discontinuities ($x_j^t\neq T_{ij}$ for all $t\in\N$,
$i=1,N$ and $j\in I(i)$). The first property is Lyapunov stability
(robustness with respect to changes in initial conditions). For
simplicity, let $\bar{x}=\{\bar{x}_i\}_{i=1,N}$ be a fixed point not
intersecting discontinuities and let
$\delta={\displaystyle\min_{i=1,N,\ j\in
    I(i)}}|\bar{x}_j-T_{ij}|>0$. Then $\bar{x}$ is asymptotically
stable and its immediate basin of attraction contains the following
cube 
\[
\{x\in\R^N\ :\ |x_i-\bar{x}_i|<\delta,\ i=1,N\}.
\]
The second property is structural stability (robustness with respect
to changes in parameters). For simplicity, assume once again that the
fixed point $\bar{x}$ exists for the parameters $(a,K_{ij},T_{ij})$
and let $\bar{\theta}$ be the corresponding code.\footnote{Every fixed
  point belongs to the attractor and the corresponding code is a
  constant sequence.} In other words, assume that $\bar{x}$ computed
by using relation (\ref{GLOBORB}) with $\bar{\theta}$, namely
\begin{equation}
\bar{x}_i=(1-a)\sum_{k=0}^{+\infty}a^k\sum_{j\in
  I(i)}K_{ij}\bar{\theta}_{ij}=\sum_{j\in
  I(i)}K_{ij}\bar{\theta}_{ij}, 
\label{FIX}
\end{equation}
satisfies the relation $\bar{\theta}_{ij}=H(s_{ij}(\bar{x}_j-T_{ij}))$
(admissibility condition (\ref{ADMICOND}) for constant codes).
If $\bar{x}$ does not intersect discontinuities
($\delta={\displaystyle\min_{i=1,N,\ j\in
    I(i)}}|\bar{x}_j-T_{ij}|>0$), then for any threshold set
$\{T'_{ij}\}$ such that $|T'_{ij}-T_{ij}|<\delta$, we have
$\bar{\theta}_{ij}=H(s_{ij}(\bar{x}_j-T'_{ij}))$ and the fixed point
$\bar{x}$ persists for the parameters $(a,K_{ij},T'_{ij})$.  Moreover,
the fixed point expression (\ref{FIX}) depends continuously on the
parameters $a$ and $K_{ij}$. Therefore, the code $\bar{\theta}$ also
satisfies the admissibility condition
$\bar{\theta}_{ij}=H(s_{ij}(\bar{x}_j-T'_{ij}))$ for the parameters
$(a',K'_{ij},T'_{ij})$ sufficiently close to $(a,K_{ij},T_{ij})$. In
short terms, every fixed point not intersecting discontinuities can be
continued for small perturbations of parameters.

Both Lyapunov stability and structural stability extend to any periodic
orbit not intersecting discontinuities, see section \ref{M-STABI} for
complete statements and further perturbation results. These properties
depend only on the distance between orbits and
discontinuities. Therefore, they may apply uniformly to all orbits in
the attractor in which case the complete dynamics is robust under small
perturbations. 

\subsection{Circuits and their symmetries}\label{R-SYMET}
\begin{figure}
\begin{center}
\input{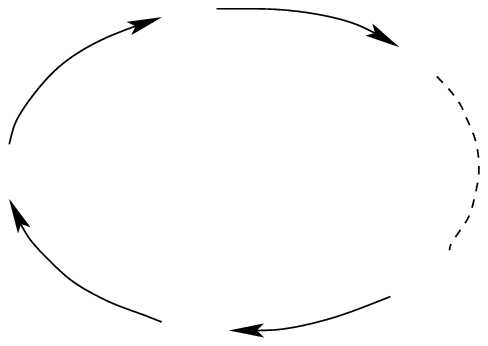}
\end{center}
\caption{A feedback circuit of $N$ genes.}  
\label{FIGCIR}
\end{figure}

The simplest regulatory networks are feedback circuits whose
graphs consist of periodic cycles of unidirectional interactions. In the
present formalism a circuit of length $N$ ($N$-circuit) can be
represented by a periodic network with $N$ genes where $I(i)=\{i-1\}$
for all $i\in\Z/N\Z$ (Figure \ref{FIGCIR}). 
When simplifying notations, in $N$-circuits the relation (\ref{MODEL}) becomes 
\begin{equation}
x_i^{t+1}=ax_i^t+(1-a)H(s_{i-1}(x_{i-1}^t-T_{i-1})),\quad i\in\Z/N\Z.
\label{CIRCUIT}
\end{equation}
Circuits are not only interesting in their own. The dynamics of
arbitrary networks can be described, in some regions of parameters, as
the combination of dynamics of independent circuits \cite{TR99}, see
section \ref{M-STABI} for detailed statements. 

A $N$-circuit is specified by interaction signs $\{s_i\}_{i\in\Z/N\Z}$, by
interaction thresholds $\{T_i\}_{i\in\Z/N\Z}$ and by the degradation rate
$a$. With as many parameters as they are, the number of cases to be
investigated is large. However, symmetry transformations apply which
allow us to considerably reduce this number.\footnote{By symmetries,
  we mean transformations acting on the original network, not on the
  subsequent dynamical graphs as in \cite{EG00}.} We present
separately flips of interaction signs and transformations of
interaction thresholds (internal symmetries). The former reduce the
number of interaction sign vectors to be considered whereas the latter
reduce the threshold domains to be studied. Details of the related
mathematical analysis are given section \ref{M-SYMET}.

\subsubsection{Flipping interaction signs}\label{R-SYMET1}
\begin{figure}
\begin{center}
\input{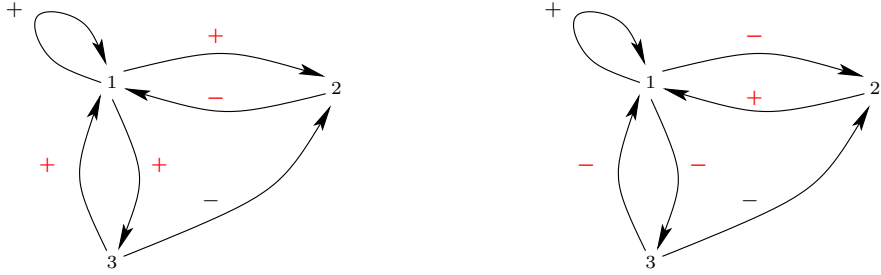}
\end{center}
\caption{A network with 3 nodes and the corresponding network obtained
  by flipping the signs of incoming and outgoing arrows from the node
  1, excepted the self-interaction sign.} 
\label{FIGFLIP}
\end{figure}

The Heaviside function possesses the following (quasi-)symmetry:
\[
H(-x)=1-H(x)\ \text{for\ all}\ x\neq 0.
\]
This property has a consequence on the dynamics of arbitrary networks
(and not only of circuits): {\em The mapping obtained by flipping the
  sign $s_{ij}$ of every incoming and every outgoing arrow from a
  fixed node, with the exception of any self-interaction, has (almost)
  the same dynamics as the original model.} Precisely, to every orbit
of the original mapping not intersecting discontinuities corresponds a
unique orbit of the new mapping (also not intersecting
discontinuities). In particular, when the attractors do not intersect
discontinuities, the asymptotic dynamics are topologically
conjugated. We refer to Lemma \ref{LemFLIP} for a complete statement
and to Figure \ref{FIGFLIP} for an example of two related networks. 

In circuits, this result implies that (with the exception of orbits
intersecting discontinuities) the dynamics only depends on the product
of signs ${\displaystyle\prod_{i\in\Z/N\Z}}s_i$ (positive or negative
circuit); a property which has been 
largely acknowledged in the literature, see e.g.\ \cite{TA90}.  

\subsubsection{Internal circuit symmetries}\label{R-SYMET2}
Lemma \ref{LemFLIP} in section \ref{M-SYMET} does not only serve to
match the dynamics of circuits with flipped signs. It can be also
applied to deduce a parameter symmetry in a circuit with fixed
signs. Let $S$ be the symmetry in $\R^N$ with respect to the point
with all coordinates equal to $\frac{1}{2}$ ($x_i=\frac{1}{2}$ for all
$i\in\Z/N\Z$). {\em The image by $S$ of a circuit orbit not 
intersecting discontinuities which exists for thresholds
$T=\{T_i\}_{i\in\Z/N\Z}$ is an orbit not intersecting discontinuities which
exists for thresholds $S(T)$.}  
 
Depending on the product of signs, other symmetries follow from
essentially cyclic permutations. The map $R$ defined by
$(Rx)_i=x_{i-1}$ is a cyclic permutation in
$\R^N$. As a representative of positive $N$-circuits we consider the
$N$-circuit with all interactions signs equal to 1: {\em The image by
  $R$ of an orbit not intersecting discontinuities which exists for 
  the thresholds $T=\{T_i\}$ is an orbit not intersecting
  discontinuities which exists for the thresholds $R(T)$}. By
repeating the argument, additional orbits $\{R^k(x^t)\}$ (unless
$R^k(x^t)=x^t$ for some $k=1,N-1$) can be obtained in this circuit. 

As a representative of negative $N$-circuits we consider a $N$-circuit
with all signs equal to 1, excepted $s_N=-1$. Let $\sigma$ be the
symmetry with respect to the hyperplane $x_1=\frac{1}{2}$: {\em The
  image by $\sigma\circ R$ of an orbit not intersecting
  discontinuities which exists for thresholds $T$ is an orbit not 
  intersecting discontinuities which exists for thresholds
  $(\sigma\circ R)(T)$}. As before, additional orbits can be obtained
in this circuit (provided that the original orbit has low or no
symmetry) by applying $\sigma\circ R$ repeatedly.

\subsection{Ghost orbits}\label{R-GHOST}
\begin{figure}
\begin{center}
\input{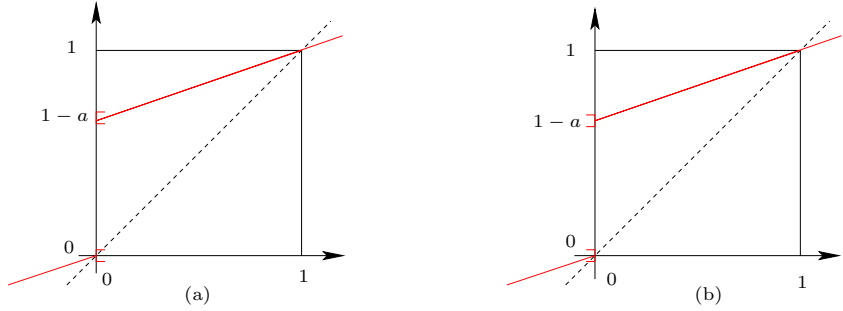}
\end{center}
\caption{(a) Graph of the map $F(x)=ax+(1-a)H(x)$. The point 0 is a
  ghost fixed point. It attracts all initial conditions $x<0$ but is
  not a fixed point. The point 1, obtained from 0 by applying the
  symmetry $S$ is a true fixed point. (b) Assuming $H(0)=0$ instead of
  $H(0)=1$ in the definition of $F$ changes the ghost fixed point to a
  genuine fixed point.} 
\label{GHOST}
\end{figure}

In the previous section, the condition of non-intersection of orbits
with discontinuities is due to lack of symmetry of the
Heaviside function at the origin ($H(-0)\neq 1-H(0)$). Applying a
symmetry transformation to an orbit intersecting some discontinuities
may result in a {\bf ghost orbit} (and vice-versa).
A ghost orbit is a sequence in phase space, which is not an orbit of
$F$, but which would be an orbit of a suitable alteration of $F$ on
some discontinuities (alterations which consists in letting $H(0)=0$
instead of $H(0)=1$). In particular, a ghost orbit must intersect some
discontinuities. The simplest example is the ghost fixed point 0 which occurs for $T=0$ in the self-activator (section \ref{R-SELFACT}), see
Figure \ref{GHOST}. (For every $T>0$, the point 0 is a stable fixed point.)

In spite of not being orbits, ghost orbits may be relevant in
applications when they attract open sets of initial conditions. This
is the case of the ghost orbits in circuits analysed below. These
are periodic sequences which exist for parameters in some boundaries
of domains where periodic orbits exist with the same symbolic
sequence. Consequently, they occur for exceptional values of
parameters. For parameters inside the domains, the periodic orbits do
not intersect discontinuities and are Lyapunov stable (see section
\ref{R-RELEV}). Ghost orbits inherit this stability with a restricted
basin of attraction (semi-neighbourhood).\footnote{For the ghost fixed point 0 in the self-activator, the restricted basin of attraction is the half-line $(-\infty,0[$, see in Figure \ref{GHOST}.}

\subsection{Comparison with other models: time delays}\label{R-DELAY}
For $a=0$, the model (\ref{MODEL}) becomes equivalent to a dynamical system on a (finite) discrete phase space, a so-called "logical network" in the literature of regulatory process models \cite{G75,T73}. Indeed for $a=0$ the symbols $\theta_{ij}^{t+1}$, can be computed by using only the symbols
$\theta_{ij}^t$ and not the variables $x_j^t$ themselves. This is because the quantities $x_i^{t+1}$ themselves can be computed by only using the symbols
$\theta_{ij}^t$.

For $a>0$, the model is no longer equivalent to a boolean network. Indeed, for $a>0$ one needs the knowledge of the variables $x_j^t$ -- and not only of the symbols -- in order to compute the next state $x_i^{t+1}$. (It may happen however that the system restricted to its attractor is equivalent to a dynamical system on a finite state space.)
\bigskip

In any case, the model (\ref{MODEL}) can be viewed as a discrete time analogue of the following system of delay differential equations
\[
\frac{dx_i}{dt}=-x_i(t)+\sum_{j\in I(i)} K_{ij}
H(s_{ij}(x_j(t-\tau)-T_{ij}))
\]
Indeed, the behaviours of (\ref{MODEL}) and of the system of delay differential equations ($\tau >0$) are remarkably similar. This is confirmed by results of the dynamics on circuits in the section below. On the opposite, in general, there are qualitative differences between the behaviours of the model (\ref{MODEL}) and of the corresponding differential equation without delay ($\tau =0$).

\noindent
More generally, the differences in behaviour between systems with and without delays and the similarities in behaviour between delay differential equations and discrete dynamical systems have been thoroughly analysed in the literature, see e.g.\ for introductory textbooks \cite{E96,KL93}. In particular, delays are known to result in the occurrence of oscillations in systems which, without delays, have only stationary asymptotic solutions. We refer to \cite{M03} for numerical example (in agreement with experimental results) of oscillations induced by delay in a negative feedback genetic regulatory circuit.
\noindent
As far as the modelisation of a real system is concerned, representing the dynamics of a regulatory network by a model with delays is relevant when the response times of various processes involved in the regulation stages (transcription, translation, etc.) are not negligible (again see \cite{M03} for a relevant example). 

\noindent
In the model (\ref{MODEL}), the degradation rate $a$ plays the role of a delay parameter: the smaller $a$, the stronger the delay is. The examples below indeed show that the oscillations (which are absent in the differential equation) are more likely to occur when $a$ is small. On the opposite, when $a$ tends to 1, the dynamics coincides with the dynamics of differential equation.

\section{Dynamics of simplest circuits}\label{R-SIMPL}
In this section, results of dynamical analysis of the four simplest
circuits are presented. For the 1-circuits and the positive 2-circuit,
the analysis is complete both in phase space and in parameter
space. Due to a rich and elaborated phenomenology, the results on the negative 2-circuit are however only partial. Still, they provide a description of the dynamics over large domains of phase space and of parameter space.

\subsection{The self-activator}\label{R-SELFACT}
The simplest network is the self-activator, namely the circuit
with one node and positive self-interaction. In this case, the mapping
$F$ becomes 
the one-dimensional map $F(x)=ax+(1-a)H(x-T)$ (where $0<T<1$) which
has very simple dynamics. Either $x<T$ and the subsequent orbit
exponentially converges to the fixed point 0. Or $x\geq T$ and the
orbit exponentially converges to the fixed point 1. 

\noindent
In terms of symbolic graph, this means that the only possible
symbolic sequences are the sequence with all symbols equal to 0
(resp.\ to 1). The first sequence corresponds to the fixed point 0, the
second to 1. Obviously, both sequences are always admissible.

\noindent
The dynamics of the self-activator is the same as those of the corresponding boolean network and of the corresponding ordinary
differential equation \cite{TA90}. In all cases, the self-activator is a bistable
system. In short terms, delays have no qualitative influence on the
self-activator. In this one-dimensional system, the reason is that, independently of the presence of delays, every orbit stays forever in its original atom and no orbit crosses the discontinuity.

\subsection{The self-inhibitor}\label{R-CONTRAC}
\begin{figure}
\begin{center}
\input{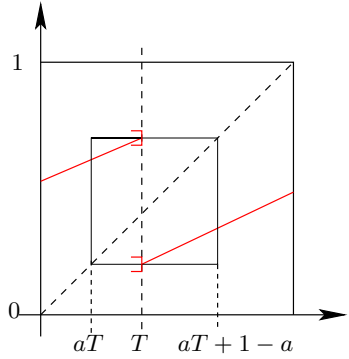}
\end{center}
\caption{Graph of the map $F(x)=ax+(1-a)H(T-x)$ together with the
  invariant absorbing interval $(aT,aT+1-a]$.} 
\label{PWCONT}
\end{figure}

For the self-inhibitor (the circuit with one node and negative
self-interaction) the mapping becomes $F(x)=ax+(1-a)H(T-x)$ whose
graph is given in Figure \ref{PWCONT}. 
Its asymptotic dynamics turns out to be the same as that of a piecewise affine
contracting rotation. Piecewise affine contracting rotations have been
thoroughly investigated \cite{BC00,C99,GT88}. The results 
presented here are immediate consequences of those in \cite{C99}. The
mathematical details are given in section \ref{M-CONTRAC}.

As indicated Figure \ref{PWCONT}, the iterations of
every initial condition in $[0,T]$ cannot stay forever in this set and
eventually enter in the complementary interval $(T,1]$ (in the sub-interval
$(T,aT+1-a]$ precisely). Conversely, every point $x^0\in (T,1]$ has an
iteration $x^t\in [0,T]$ (which indeed belongs to the
sub-interval $(aT,T]$). Therefore all orbits oscillate between
the two (sub-)intervals. 

\noindent
In terms of symbolic dynamics, it means that the symbolic graph is the
complete graph (that is to say 0 and 1 can both be followed by either 0
or 1). However, only special paths correspond to admissible
sequences. Indeed, the analysis of the admissibility 
condition in this case proves that the only admissible sequences are
the codes generated by a rigid rotation on the circle
($x\mapsto x+\nu\ \text{mod}\ 1$) with a unique rotation number $\nu$. 

\noindent
Back to the phase space, the corresponding orbits themselves are
given, up to a change of variable, by such a rigid rotation. Moreover,
they attract all initial conditions (Theorem \ref{CONJUG}). Therefore,
the asymptotic 
dynamics is entirely characterised by the rotation number which corresponds to the mean
fraction of iterations spent in the interval $(T,aT+1-a]$.

\noindent
For almost all values of parameters, the asymptotic orbits are 
genuine orbits. But in a set of parameters $(a,T)$ with zero
Lebesgue measure, all orbits approach a unique ghost periodic
orbit.\footnote{In this case, the rotation number is still well-defined - and
is a rational number - but the attractor is empty.}

\begin{figure}
\begin{center}
\input{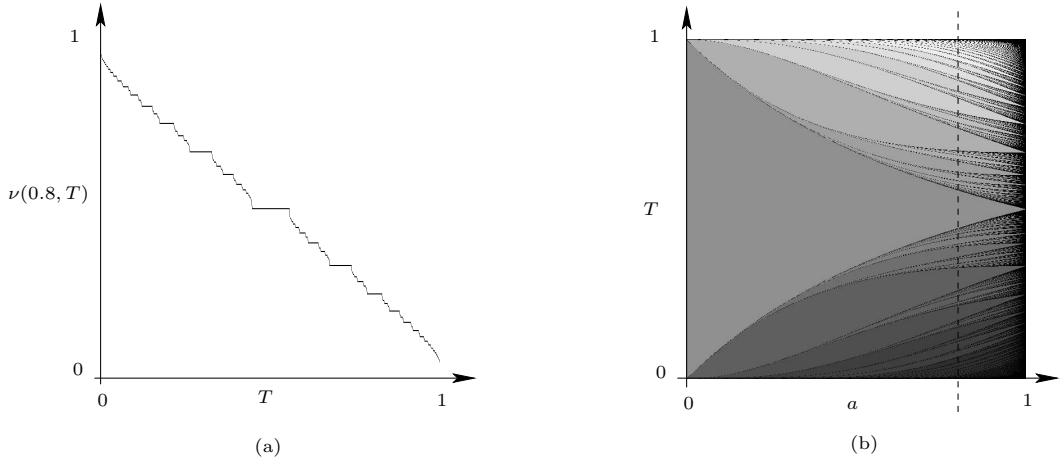}
\end{center}
\caption{The self-inhibitor. (a) Graph of the map $T\mapsto \nu(a,T)$
  for $a=0.8$. (b) Gray-level plot of the map $(a,T)\mapsto\nu(a,T)$
  (White = 0 - Black = 1).} 
\label{ROTANUMB}
\end{figure}

The rotation number $\nu(a,T)$ depends continuously on $a$ and on $T$
(small changes in parameters induce small changes in the rotation
number). Moreover for $a>0$ all maps $T\mapsto \nu(a,T)$ are decreasing with
range $(0,1)$ and have a peculiar structure called a Devil's 
staircase, see Figure \ref{ROTANUMB} (a).

\noindent
This structure combines continuity and the
existence of intervals (plateaus) where the rotation number is a
constant rational number. Plateaus are due to structural 
stability of orbits not intersecting discontinuities (section
\ref{R-RELEV}). Indeed, when the rotation 
number is rational, the map $F$ has a periodic orbit which persists
under small (suitable) perturbations of parameters. Since the rotation
number does not depend on the initial condition, it remains constant
while this periodic orbit persists and we have a
plateau.\footnote{Ghost periodic orbits occur for $T$ at the right
  boundary of plateaus.} 

\noindent
Additional properties of the rotation number are the symmetry
$\nu(a,1-T)=1-\nu(a,T)$ which is a consequence of the symmetry map $S$
(section \ref{R-SYMET2}) and the unique plateau $\nu(0,T)=\frac{1}{2}$
for $a=0$ (Figure \ref{ROTANUMB} (b)). Therefore, for $a=0$ and every
$T\in (0,1)$, every orbit asymptotically approaches a unique
2-periodic orbit and the attractor is just the same as in the
corresponding boolean model \cite{TA90}. 

\begin{figure}
\begin{center}
\input{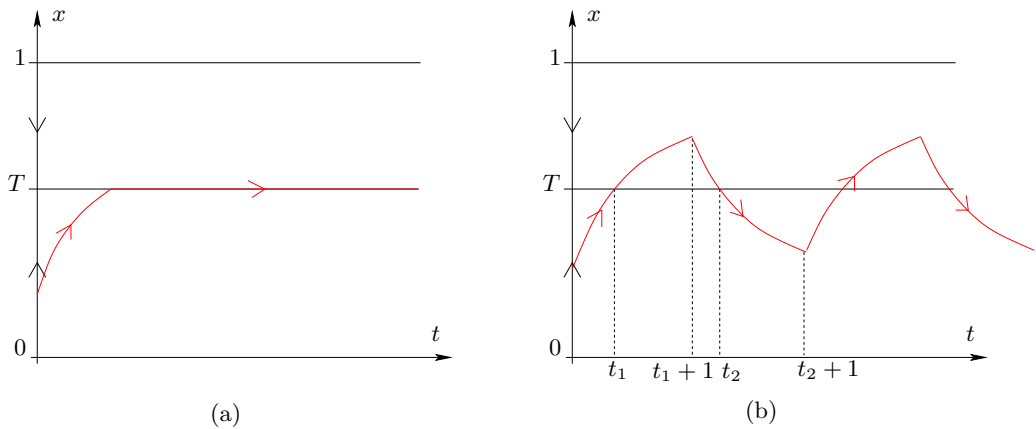}
\end{center}
\caption{Trajectories in the continuous time self-inhibitor without
  delay (a) and with delay $\tau =1$ (b).} 
\label{DELAY}
\end{figure}
As argued in section \ref{R-DELAY}, the permanent oscillations in the
dynamics of $F$ can be attributed to the presence of delays. This is
confirmed by the oscillating behaviour of solutions of the
delay differential equation 
\[
\frac{dx}{dt}=-x(t)+H(T-x(t-1))
\]
see Figure \ref{DELAY} (b). The behaviour of the corresponding
system without delay is qualitatively different. Indeed for any $T\geq
0$ every trajectory of the ordinary differential equation 
\[
\frac{dx}{dt}=-x(t)+H(T-x(t))
\]
converges to the globally attracting stationary point $x=T$ (see
Figure \ref{DELAY} (a)).\footnote{Defining $H(0)=T$ instead of $H(0)=1$.}

In addition, the oscillations of $F$ belong to the absorbing interval
$(aT,aT+1-a]$ and this interval reduces to the point $T$ in the limit
$a\to 1$. Under this point of view, the attractor of $F$ converges in
the limit of vanishing delay to the attractor of the differential
equation without delay.
  
\subsection{The positive 2-circuit}
\begin{figure}
\begin{center}
\input{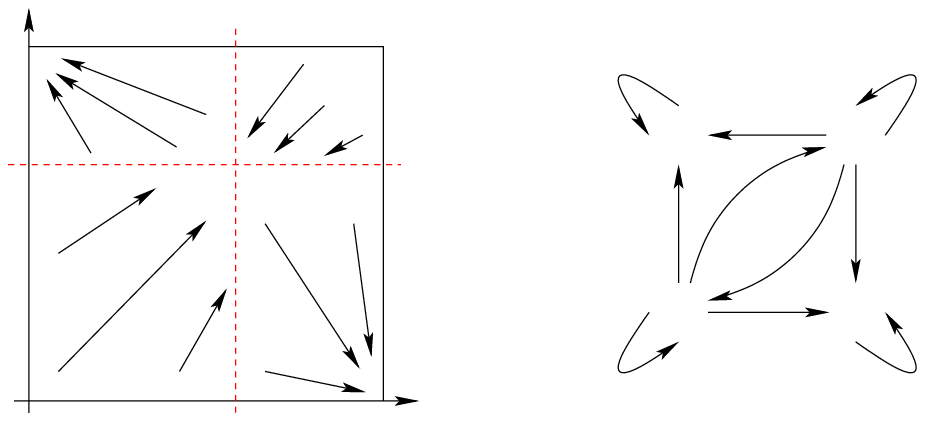}
\end{center}
\caption{The positive 2-circuit. (a) Atoms in phase space with affine
  dynamics and directions of motions. (b) Associated symbolic graph.} 
\label{POSCIR}
\end{figure}

According to the symmetry of flipping all interaction signs (section
\ref{R-SYMET1}) the positive 2-circuit can be obtained by choosing, in relation (\ref{CIRCUIT}) with $N=2$, 
either $s_1=s_2=1$ or $s_1=s_2=-1$. In order to compare with results 
on other models in the literature \cite{EG00,TA90}, we have opted for
the second choice. The system (\ref{CIRCUIT}) then becomes the following
system of cross-inhibitions
\[
\left\{\begin{array}{l}
x_1^{t+1}=ax_1^t+(1-a)H(T_2-x_2^t)\\
x_2^{t+1}=ax_2^t+(1-a)H(T_1-x_1^t)
\end{array}\right.
\]
which can be viewed as the iterations of a mapping
$F$ of the square $[0,1]^2$. For such mapping, the symbolic coding
follows from the partition of the square into 4 atoms labelled
by 00, 01, 10 and by 11, see Figure \ref{POSCIR} (a).
The corresponding symbolic graph is given Figure \ref{POSCIR}
(b). According to this graph, for every initial condition in $[0,1]^2$,
one of the following assertions holds 
\begin{itemize}
\item{} either the orbit visits only the atoms 00 and 11, 
\item{} or the orbit enters one of the atoms 01 or 10 and is trapped
  inside it forever. 
\end{itemize}
When trapped in 01 the orbit converges to $(1,0)$
(lower right corner of the square). When in 10 the orbit converges to
$(0,1)$. The fixed points $(1,0)$ and $(0,1)$ exist and are stable for any
parameter $a\in [0,1)$ and $T_1,T_2\in (0,1)$.

\begin{figure}
\begin{center}
\input{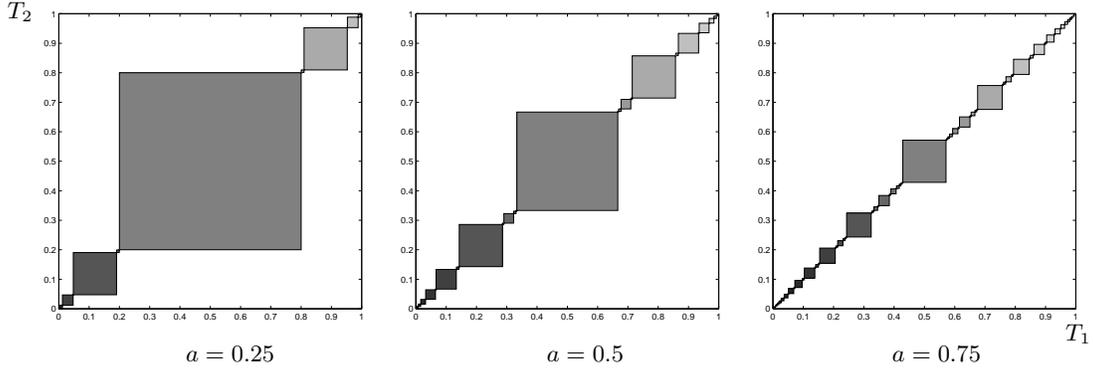}
\end{center}
\caption{The positive 2-circuit. Gray level plots of the domains in the threshold plane $(T_1,T_2)$ where the
  diagonal contains an orbit with rotation number $\nu$ (White: $\nu
  =0$ - Black: $\nu=1$). A projection of the domains in the
  plane $(a,T_1)$ (or in the plane $(a,T_2)$) is given Figure
  \ref{ROTANUMB} (b).} 
\label{SQUAR_POS}
\end{figure}

That an orbit visits only 00 and 11 depends on the parameters. If this
happens, the orbit must oscillate between the two atoms and must
asymptotically approach the diagonal, i.e.\ $x_1^t-x_2^t\to 0$ when
$t\to\infty$. \footnote{This follows from the fact that if $(x_1^t,x_2^t)\in
  00\cup 11$, then $x^{t+1}_1-x_2^{t+1}=a(x_1^t-x_2^t)$.} An orbit on
the diagonal is characterised by 
$x_1^t=x_2^t\equiv x^t$. According to the definition of $F$,
the quantity $x^t$ must satisfy the dynamics of the self-inhibitor
simultaneously for $T=T_1$ and for $T=T_2$ 
\[
x^{t+1}=ax^t+(1-a)H(T_1-x^t)=ax^t+(1-a)H(T_2-x^t)\quad t\in\N
\]
The results on the self-inhibitor imply that these equalities hold for some
$x^0\in [0,1]$ iff the thresholds $T_1$ and $T_2$ are such that
$\nu(a,T_1)=\nu(a,T_2)$. 
According to relation (\ref{INETHRE}) in section \ref{M-CONTRAC}, the attractor
(or a ghost periodic orbit) of the system of mutual inhibitions
intersects the diagonal iff $T_1$ and $T_2$ belong to the 
interval (the point if $\nu$ is irrational)
$[\overline{T}(a,\nu),\overline{T}(a,\nu-0)]$ for some $\nu\in
(0,1)$, see Figure \ref{SQUAR_POS}. 

\noindent
When a periodic orbit on the diagonal exists, its basin of attraction can be determined. For instance, there exists a domain of parameters -- defined by the inequality $a<\min\{T_1,1-T_1,T_2,1-T_2\}$ -- for which the image of 00 is contained in 11 and the image of 11 is contained in 00. In other words, for such parameters, there exists a 2-periodic orbit on the diagonal which attracts every initial condition in the atoms 00 and 11.

\noindent
As a consequence of the continuous dependence of $\nu(a,T)$ with $a$,
the squares $[\overline{T}(a,\nu),\overline{T}(a,\nu-0)]^2$ vary
continuously with $a$. For $a=0$, the central square corresponding to
$\nu=1/2$ coincides with $[0,1]^2$. Just as in the corresponding
boolean network \cite{TA90}, the positive 
2-circuit possesses for $a=0$ two stable fixed points and a stable
2-periodic orbit.  

These permanent oscillations do not occur in the corresponding
system of coupled ordinary differential equation which only present
two stable fixed points (and possibly a hyperbolic fixed point on a
separatrix) \cite{EG00,TA90}. As before, oscillations can be
attributed to a delay effect due to discreteness of time in the model
(\ref{MODEL}). That a time delay is necessary to obtain permanent
oscillations in positive circuits has already been
acknowledged in the literature \cite{G97,NT97}. 

In the limit $a\to 1$, the union of rectangles $[aT_1,T_1]\times
[aT_2,T_2]\cup [T_1,aT_1+1-a]\times [T_2,aT_2+1-a]$ which contain the
oscillations reduces to the point $(T_1,T_2)$. Once again in the
limit of vanishing delay, the attractor reduces to that of the differential
equation without delay. 

\subsection{The negative 2-circuit}\label{R-NEGAT}
\begin{figure}
\begin{center}
\input{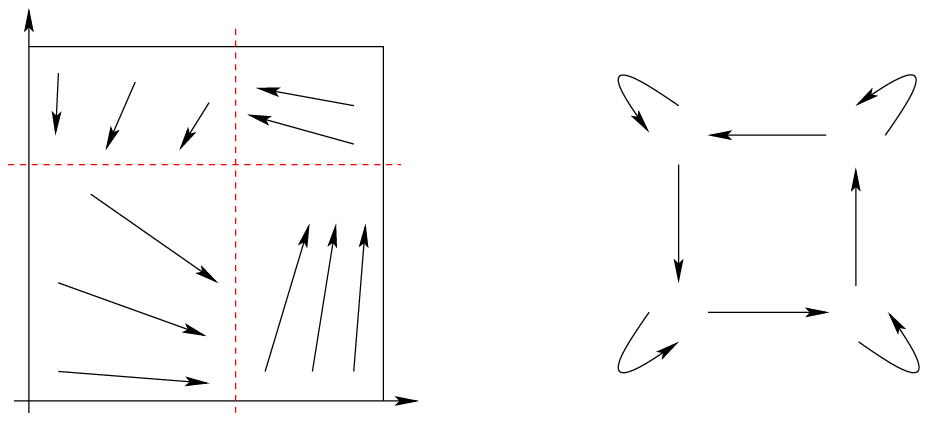}
\end{center}
\caption{The negative 2-circuit. (a) Atoms in phase space with affine
  dynamics and directions of motions. (b) Associated symbolic graph.} 
\label{NEGCIR}
\end{figure}

Up to a flip of all interaction signs, the negative 2-circuit can be
obtained by choosing $s_1=1$ and $s_2=-1$ in relation (\ref{CIRCUIT}) with
$N=2$. As for the positive circuit, coding follows from a partition of
the unit square into 4 atoms 00, 01, 10 and 11, see Figure
\ref{NEGCIR} (a). \footnote{Note however that the correspondence
  between labels and atoms differs from that in the positive circuit -
  compare Figures \ref{POSCIR} (a) and \ref{NEGCIR} (a).} The
associated symbolic graph is given Figure \ref{NEGCIR} (b). As this figure suggests, no orbit can stay forever in an arbitrary given atom and every orbit visits sequentially every atom.

\noindent
Numerical simulations indicate that {\em in most cases} of initial conditions and of parameters, this recurrence is regular, i.e.\ the orbit winds regularly around the intersection $(T_1,T_2)$ of interaction thresholds. Motivated by these numerical results, we have characterised such regular orbits and we have accomplished the mathematical analysis of their existence and of their parameter dependence.  This analysis is reported in section \ref{M-NEGAT} and its main results are presented in the next two sections.  

\subsubsection{Balanced periodic orbits}
According to the symbolic graph, the simplest regular behaviour is a periodic
orbit passing the same number $p$ of consecutive steps in
each atom ({\bf balanced orbit}). Its code is given by $(00^p\ 01^p\ 11^p\ 10^p)^\infty$ or formally \footnote{The
  symbol $\theta^t$ denotes the pair $(\theta_1^t,\theta_2^t)$. The
  map $\sigma\circ R$ defined in section \ref{R-SYMET2} simply becomes
  the rotation by $\frac{\pi}{2}$ with centre
  $(\frac{1}{2},\frac{1}{2})$ in the present case. It
  transforms $(x_1,x_2)\in [0,1]^2$ into $(1-x_2,x_1)$. As a consequence, we have $S=(\sigma\circ R)^2$, i.e.\ there is indeed only one (independent) internal symmetry in the negative 2-circuit.}
\[
\theta^{t+p}=(\sigma\circ R)(\theta^t)\quad \text{for all}\
t\in\Z\quad\text{and}\quad \theta^t=00\quad \text{for all}\ t=1,\cdots, p
\]
According to expression (\ref{GLOBORB}) section \ref{M-ATTRA}, the 
corresponding orbit has the same symmetry. Namely, for all $t\in\Z$,
$x^{t+p}$ is the image of $x^t$ under the rotation by angle
$\frac{\pi}{2}$ with centre $(\frac{1}{2},\frac{1}{2})$. In
particular the orbit is $4p$-periodic and its
components $x^t$ lie on the boundaries of a square.
\begin{figure}
\begin{center}
\input{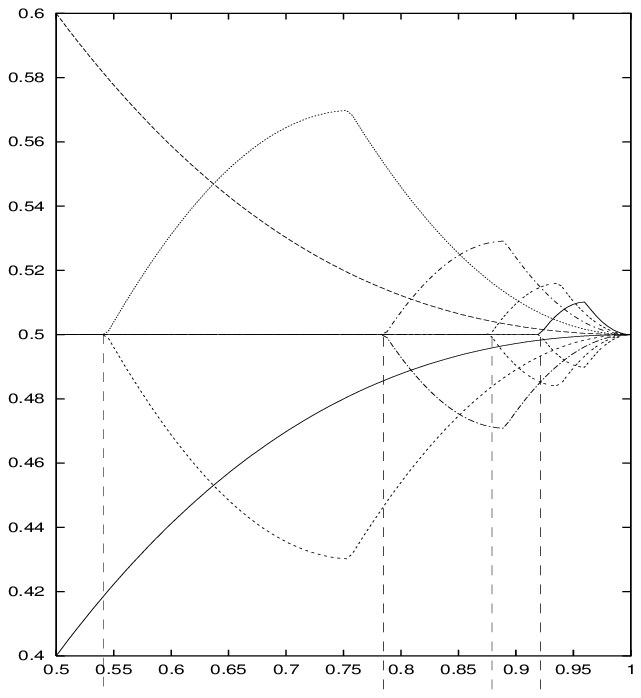}
\end{center}
\caption{The negative 2-circuit. Projections of parameter domains of
  existence of balanced orbits -- for $p=1,5$ -- on the plane
  $(a,T_1)$ (or by symmetry on the plane $(a,T_2)$). The
  representation here differs from Figure \ref{SQUAR_POS}
  (intersections with planes $a=\text{constant}$ in threshold
  space). A representation as in Figure \ref{SQUAR_POS} would have 
  presented nested squares centred at $(\frac{1}{2},\frac{1}{2})$ and
  symmetric with respect to the diagonal; the number and the size of
  squares depending on $a$.}  
\label{PARADOM}
\end{figure}

\noindent
The existence domain of an arbitrary balanced orbit have been computed explicitely, see third item in section \ref{EXAMPL-NEGAT}. These domain have been represented on Figure \ref{PARADOM}.

\noindent
The product structure in threshold space and the rotation 
symmetry $\sigma\circ R$ imply that, when non-empty, the existence domain in the 
plane $(T_1,T_2)$ is a square centred at $(\frac{1}{2},\frac{1}{2})$,
symmetric with respect to the diagonal and which depends on $a$. 
For $p=1$ the square is non-empty for every $a\in [0,1)$ and fills the whole
square $(0,1)^2$ for $a=0$. For every $p>1$, the square is non-empty
iff $a>a_p$. The critical value $a_p>\frac{1}{2}$, increases with $p$ and converges to 1 when $p\to\infty$. 

\begin{figure}
\epsfxsize=7truecm 
\centerline{\epsfbox{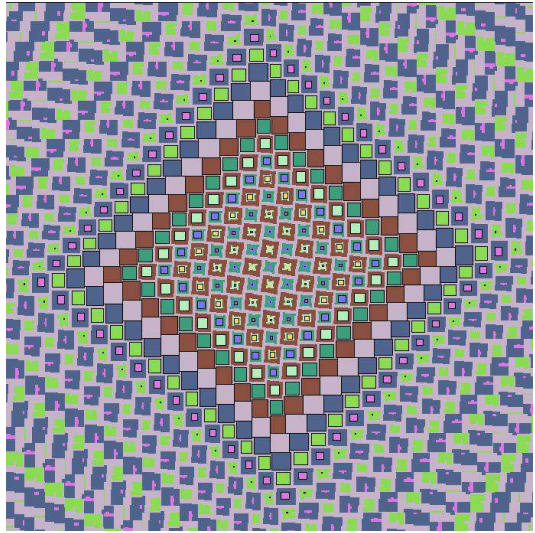}}
\caption{Colour plot of basins of attractions in the square $[0.425,0.575]^2$ for the negative 2-circuit for $a=0.99$ and $T_1=T_2=\frac{1}{2}$ -- obtained from a numerical simulation. The picture shows that there are only balanced orbits. Several immediate basins -- characterised by products of intervals -- clearly appear. For instance, sea-green squares correspond to the balanced orbit with $p=8$, brown squares to the orbit with $p=9$, thistle squares to the orbit with $p=10$, etc. In the present case, the balanced $4p$-periodic orbits are known to exist for $p=1,\cdots,14$ (indeed we have $a_{14}<0.99$).}
\label{BASIN}
\end{figure}

\begin{figure}
\epsfxsize=7truecm 
\centerline{\epsfbox{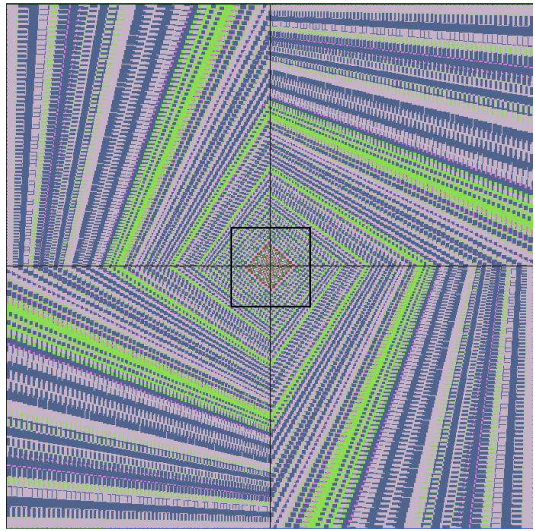}}
\caption{Colour plot of the basins of attraction in the square $[0,1]^2$ (complete phase space) for the negative 2-circuit with parameters as in Figure \ref{BASIN} (i.e.\ $a=0.99$ and $T_1=T_2=\frac{1}{2}$). The picture clearly shows that any orbit has a spiral trajectory toward some balanced orbit located in the region -- delimited by the square -- displayed in Figure \ref{BASIN}.}
\label{SPIRAL}
\end{figure}

\noindent
For $a=0$, the 4-periodic orbit is the unique orbit in the
attractor for all pairs $(T_1,T_2)$. The dynamics is equivalent to the
corresponding boolean network \cite{TA90}. On the opposite, for $a$ 
arbitrarily close to 1 and $T_1=T_2=\frac{1}{2}$, an arbitrary large 
number of stable balanced orbits coexist and we have multi-stability, see
Figure \ref{BASIN}. 
Their components tend to $(\frac{1}{2},\frac{1}{2})$ when $a$ tends to 1. Thus, when $a$ is close to 1 and $T_1=T_2=\frac{1}{2}$, at large scales, the dynamics is as for the corresponding system of coupled differential equations (see \cite{TA90} for an analysis of such system). The attractor is (concentrated in a neighbourhood of) the point
$(\frac{1}{2},\frac{1}{2})$ and every orbit has a spiral trajectory toward this point (region), see Figure \ref{SPIRAL}. Multi-stability occurs at a smaller scale, inside this neighbourhood of $(\frac{1}{2},\frac{1}{2})$.

\subsubsection{Regular orbits}
\begin{figure}
\begin{center}
\input{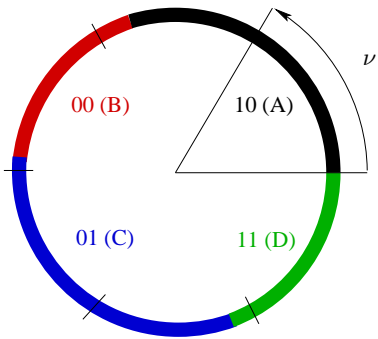}
\end{center}
\caption{A regular orbit code is generated by a rigid rotation
on the unit circle $x\mapsto x+\nu\text{ mod }1$ composed by 4
arcs. Each arc (specified by a colour) is associated with an 
atom of the partition. Here $\nu=\frac{1}{6}$ and the generated code is
$(10^2\ 00\ 01^2\ 11)^\infty$ and is 6-periodic.} 
\label{REGULAR}
\end{figure}

Balanced periodic orbits are special cases of orbits winding
with a regular motion around the point $(T_1,T_2)$. An orbit is said to wind regularly around $(T_1,T_2)$ ({\bf regular orbit}) if its code is generated by the
orbit of the rotation $x\mapsto x+\nu\text{ mod }1$ on the unit circle
composed by 4 arcs; each arc being associated with an atom 00, 01, 10 or 11,
see Figure \ref{REGULAR}. 

\noindent
In short terms the code of a regular orbit is characterised by the 5-uple 
$(A,B,C,D,\nu)$ where $A$ (resp.\ $B$, $C$, $D$) is the length of the arc
associated with 10 (resp.\ 00, 01, 11) and $\nu$ is the rotation
number. In particular the balanced $4p$-periodic orbit is a regular orbit
for which the 4 arc lengths are equal to $\frac{1}{4}$ and the rotation number is 
equal to $\frac{1}{4p}$.

\noindent
The balanced $4p$-periodic orbits exist when both
thresholds $T_1,T_2$ are sufficiently close to $\frac{1}{2}$ -- see Figure \ref{PARADOM}. When this is not the case, the orbits spend more iterations in some atom(s) than in others. The simplest case is when the orbit spends the same number of steps per winding in 3 atoms and a different number of steps per winding in the fourth atom (3 arc lengths are equal). 

\noindent
In order to cover larger parameter domains, we consider those regular orbits with only 2 (consecutive) equal arc lengths, say $A=B$. Moreover, these lengths and the rotation number are chosen so that the orbits spend one step per winding in each the corresponding atoms (10 and 00). 

\noindent
Motivated by changes in dynamics with parameters, instead of considering regular orbits with "isolated" numbers of steps per winding in the two remaining atoms, we consider the families of regular orbits -- called {\bf $(p,\rho)$-regular orbits} -- for which one of these numbers is a fixed arbitrary integer and for which the other number is a continuous parameter.\footnote{When this number is not an integer, it is to be interpreted as a mean number of steps spent per winding in a given atom.} That is to say the length $C$ and the rotation number are chosen so that the number of iterations spent per winding in $01$ is $p\in\N$. The (average) number of iterations spent per winding in $11$ is represented by the real number $\rho\geq 1$. \footnote{Technically speaking, we have $A=B=\nu$, $C=p\nu$ and by normalisation $D=1-(p+2)\nu$ which implies that $\rho:=\frac{D}{\nu}=\frac{1}{\nu}-(p+2)$. We refer to section \ref{M-NEGAT} for more details.}
\bigskip

By analysing the admissibility condition of the corresponding codes, we have obtained
the following results on the existence, uniqueness and parameter dependence of
$(p,\rho)$-orbits. (The proof is given in section \ref{EXAMPL-NEGAT}, second item.) We need the unique real root of the polynomial $a^3+a^2+a-1$ -- denoted by $a_c$ -- which is positive ($a_c\sim 0.544$). 
\begin{Thm} {\bf (Families of regular orbits and their parameter dependence. Simple case)}
The $(p,\rho)$-regular orbit exists iff $(T_1,T_2)$ belongs to a unique rectangle $I_1(a,p,\rho)\times I_2(a,p,\rho)$ which exists for every $p\geq 1$ and
$\rho\geq 1$ provided that $a\in (0,a_c]$. 

\noindent
The boundaries of the intervals $I_i(a,p,\rho)$ $(i=1,2)$ are strictly increasing functions of $\rho$. The boundaries of $I_2(a,p,\rho)$ tend to 1 when $\rho$ tends to $\infty$. Moreover $I_2(a,p,\rho)$ reduces to a point iff $\rho$ is irrational.

\noindent
The intervals $I_1(a,p,\rho)$ and $I_1(a,p,\rho')$ intersect when $\rho$ and $\rho'$are sufficiently close. On the other hand, we have $I_2(a,p,\rho)<I_2(a,p,\rho')$ whenever $\rho<\rho'$ and the union ${\displaystyle\bigcup_{\rho\geq 1}}I_2(a,n,\alpha)$ consists of an interval excepted a countable nowhere dense set (where we have a ghost regular periodic orbit instead).  
\label{STATEREG}
\end{Thm}
In other words, when $\rho'>\rho$, the existence domain in the threshold plane of the $(p,\rho)$-regular orbit (the rectangle $I_1(a,p,\rho')\times I_2(a,p,\rho')$) lies (strictly) above and at the right of the existence domain of the $(p,\rho)$-regular orbit (the rectangle $I_1(a,p,\rho)\times I_2(a,p,\rho)$).
As a consequence given $(a,p,T_1,T_2)$ the number $\rho$ is unique. Moreover it is an increasing function of $T_2$, with a Devil's staircase structure, and which tends to $\infty$ when $T_2$ tends to 1. 
\bigskip

\begin{figure}
\begin{center}
\input{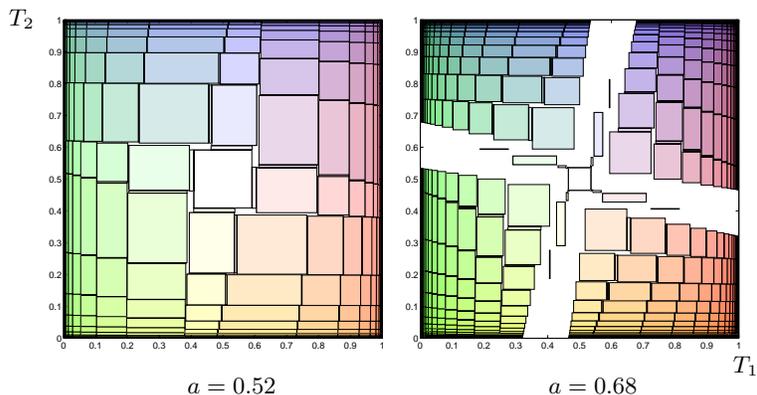}
\end{center}
\caption{Colour plots of the existence domains in the threshold plane $(T_1,T_2)$ of some regular orbits for $a=0.52$ and $a=0.68$. The orbits are characterised by one iteration per winding in any two consecutive atoms, $p$ iteration(s) per winding in a third atom and $\rho$ iteration(s) per winding in the remaining atom. On each picture, the domains where all orbits pass one iteration in $10$ and in $00$ -- the $(p,\rho)$-regular orbit -- concern the right upper quadrant (red and blue). Light (resp.\ dark) red corresponds to small (resp.\ large) $p$, i.e.\ number of iteration(s) per winding in $01$. Light (resp.\ dark) blue corresponds to small (resp.\ large) $\rho$, i.e.\ (mean) number of iteration(s) per winding in $11$. The domains where the orbits pass one iteration in $00$ and in $01$ concern the left up quadrant, etc. See text for further details.}
\label{SQUAR_NEG}
\end{figure}

The expression of the rectangle boundaries are explicitely known (see section \ref{M-NEGAT}). By numerically computing these quantities and by applying the symmetry $\sigma\circ R$, we have obtained Figure \ref{SQUAR_NEG}. This figure presents the existence domains in threshold space, for $a=0.52$ and for $a=0.68$ of all regular orbits which pass one iteration per winding in any two consecutive atoms, $p$ iteration(s) per winding in a third atom and $\rho$ iteration per winding in the remaining atom.

\noindent
On the first picture of Figure \ref{SQUAR_NEG} ($a=0.52<a_c$), the central square (in white) is the existence domain of the balanced 4-periodic orbit. The series of rectangles above it (from light to dark blue) are existence domains of $(1,\rho)$-regular orbits. In particular, the first large rectangle above the white square is the $(1,2)$-regular orbit existence domain. The small rectangle in between corresponds to the $(1,\frac{3}{2})$-regular orbit (the rectangles in between and corresponding to other $(1,\rho)$-regular orbits with $1<\rho<2$ are too thin to appear on the picture). The second large square above corresponds to the $(1,3)$-regular orbit, etc. 

\noindent
The series of rectangles extending upward at the right of the series corresponding to $(1,\rho)$-regular orbits, are existence domains of $(2,\rho)$-regular orbits. In particular, the rectangle in light red at the right of the white square corresponds to the $(2,1)$-regular orbit. Obviously, it is symmetric to the rectangle corresponding to the $(1,2)$-regular orbit. The series extending upward at the right of the series corresponding to $(2,\rho)$-regular orbits corresponds to $(3,\rho)$-regular orbits, and so on. Series are visible up to $p=7$.
\bigskip

The second picture of Figure \ref{SQUAR_NEG} ($a=0.68>a_c$) shows that when $a$ increases beyond $a_c$, some rectangle persist whereas other do not. This is confirmed by the next statement which claims that the rectangle corresponding to the pair $(p,\rho)$, with any $\rho$ sufficiently large, persist if $p$ is large and do not persist if $p$ is small. 
\begin{Pro}
Provided that $p$ is sufficiently large (i.e. for any $p$ larger than a critical value $p_a$ which depends on $a\in (0,1)$), the results of Theorem \ref{STATEREG} extend to the $(p,\rho)$-regular orbits with arbitrary $\rho$ larger than a critical value (which depends on $a$ and on $p$). 

\noindent
(Theorem \ref{STATEREG} states that $p_a=1$ and that the critical value of $\rho$ equals 1 whenever $a<a_c$.) When $a>a_c$, both $p_a$ and the critical value of $\rho$ are larger than 1. The critical integer $p_a$ tends to $\infty$ when $a$ tends to 1. 
\label{OPTIMAL}
\end{Pro}
In particular, for $a=0.68$ (second picture of Figure \ref{SQUAR_NEG}), we have $p_a=2$ and for any $p\geq 2$ the critical value of $\rho$ is (at most) 2. That is to say the $(p,\rho)$-regular orbit exists for any $p\geq 2$ and any $\rho\geq 2$. On the other hand we have $a=0.68>\overline{a}_{1,1}$ (see Theorem \ref{STATEREG2} below for the definition of $\overline{a}_{1,1}$) and thus $(p,\rho)$-regular orbit with $p=1$ only exist for isolated values of $\rho$. In particular, the balanced 4-periodic orbit exists. 
\bigskip
 
\begin{figure}
\begin{center}
\input{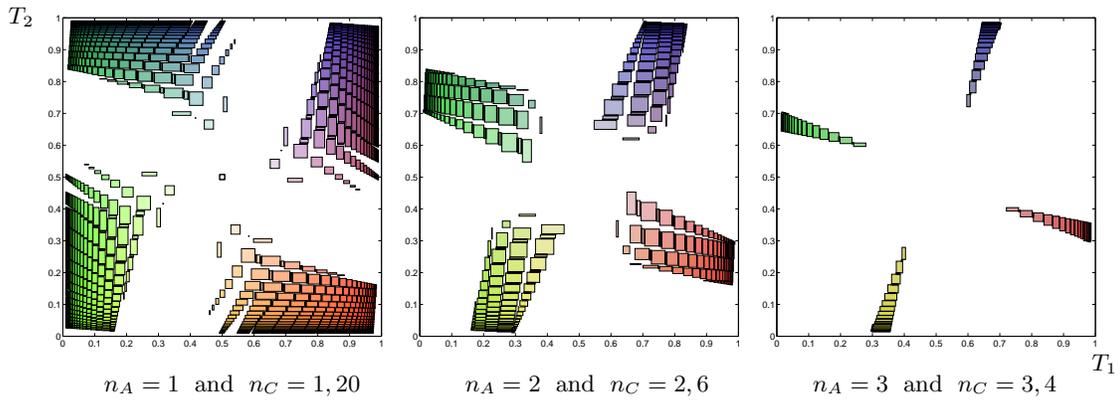}
\end{center}
\caption{Colour plots of $(n_A,1,n_C,\rho)$-regular orbit existence domains in the threshold plane for $a=0.842$. See text for details.}
\label{SQUAR_NEG2}
\end{figure}

We also have obtained results on more general regular orbits than only those spending one iteration per winding in two consecutive atoms. Inspired by the previous statement, one may wonder about the existence of regular orbits, denoted $(n_A,n_B,n_C,\rho)$-regular orbits, with an arbitrary integer number of iterations spent per winding in each of 3 atoms (say $n_A$ iterations in 10, $n_B$ iterations in $00$ and $n_C$ iterations in 01) and an arbitrary mean number of iterations spent per winding in the fourth atom ($\rho\in\R,\ \rho\geq 1$ iterations on average in 11).\footnote{Technically speaking, we have $A=n_A\nu$, $B=n_B\nu$, $C=n_C\nu$ and $\rho:=\frac{D}{\nu}=\frac{1}{\nu}-(n_A+n_B+n_C)$.} (Under this notation, the previous $(p,\rho)$-regular orbits are $(1,1,p,\rho)$-regular orbits). As claimed in the next statement, it turns out that the number $\rho$ can be arbitrary large only if the number of iterations in the opposite atom is equal to 1.
\begin{Thm} {\bf (Families of regular orbits and their parameter dependence. General case)}
Let $n_A\geq 1$ and $n_C\geq 1$ be arbitrary integers. 

\noindent
The $(n_A,n_B,n_C,\rho)$-regular orbit can exist -- upon a suitable choice of the parameters $(a,T_1,T_2)$ -- for any $\rho$ in an interval of the form $(\rho_c,\infty)$ only if $n_B=1$. 

\noindent
The $(n_A,1,n_C,\rho)$-regular orbit exists iff $(T_1,T_2)$ belongs to a unique rectangle $I_1(a,n_A,n_C,\rho)\times I_2(a,n_A,n_C,\rho)$ which exists provided that $a\in [\underline{a}_{n_A,n_C},\overline{a}_{n_A,n_C}]$ and that $\rho$ is larger than a critical value (say $\rho\geq\rho_{a,n_A,n_C}$). The numbers $\underline{a}_{n_A,n_C}$ and $\overline{a}_{n_A,n_C}$ are known explicitely. 

\noindent
The dependence of intervals $I_i(a,n_A,n_C,\rho)$ $(i=1,2)$ on $\rho$ is just as in Theorem \ref{STATEREG}.
\label{STATEREG2}
\end{Thm} 
For the proof and for explicit expressions, see section \ref{EXAMPL-NEGAT}, first item. Results of the numerical computation of $(n_A,1,n_C,\rho)$-regular orbit (and their symmetric) existence domains for $a=0.842$, for 3 values of $n_A$ and for the values of $n_C$ such that $a\in [\underline{a}_{n_A,n_C},\overline{a}_{n_A,n_C}]$ are presented on Figure \ref{SQUAR_NEG2}. On these pictures the $(n_A,1,n_C,\rho)$-regular orbit existence domains are in the right upper quadrants. The other domains corresponds to orbits obtained by applying the symmetry $\sigma\circ R$.

\noindent
On the first picture of Figure \ref{SQUAR_NEG2},  we have represented domain for $n_A=1$ and for $n_C$ running from 1 to 20. Explicit calculations show that $a=0.842>\overline{a}_{1,n_C}$ for $n_C=1,2$ and 3. No continuum of domains  but only isolated domains exist for $n_C=1$ and 2. For $n_C=3$ we have a continuous series of domains for $\rho$ inside a finite interval (first series extending upward in the right upper quadrant). On the other hand we have $a=0.842\in [\underline{a}_{1,n_C},\overline{a}_{1,n_C}]$ for all $n_C\geq 4$. Thus for such $n_C$, we have a continuous series of domains for all $\rho$ larger than a critical value (unbounded intervals). Each series for $n_C+1$ stands at the right of the series for $n_C$.

\noindent
On the second picture of Figure \ref{SQUAR_NEG2} we have $n_A=2$ and $n_C$ runs from 2 to 6. The parameter $a=0.842\in [\underline{a}_{2,n_C},\overline{a}_{2,n_C}]$ for $n_C=2,3,4,5$ and 6. The continuous series of domains for $n_C=6$ is quite small. 

\noindent
On the third picture of Figure \ref{SQUAR_NEG2} we have $n_A=3$ and $n_C=3,4$. The parameter $a=0.842\in [\underline{a}_{3,n_C},\overline{a}_{3,n_C}]$ for $n_C=3$ and 4. The continuous series for $n_C=4$ is barely visible.
\bigskip

\noindent
{\bf Families of regular orbits with arbitrary iterations per winding in the four atoms}

\noindent
The existence and the parameter dependence of arbitrary families of $(n_A,n_B,n_C,\rho)$-regular orbits ($n_A>1,n_B>1$ and $n_C>1$) where $\rho$ is varying in an interval $(\rho_1,\rho_2)$ (with $\rho_2$ necessarily finite) remain to be investigated. According to the analysis developped in section \ref{M-NEGAT}, the corresponding existence conditions are explicitely known (see Proposition \ref{FINALCOND}). Moreover the $\rho$-dependence on threshold parameters is as in Theorems \ref{STATEREG} and \ref{STATEREG2} (see Proposition \ref{ENCORE}). Thus, given $n_A,n_B,n_C,\rho_1$ and $\rho_2$, one only has to check the values of $a$ for which the existence conditions hold.  

\noindent
In the particular case where $n_A=n_B=n_C=p$, we checked the conditions numerically for several values of $p$. A statement similar to Theorem \ref{STATEREG2} results: {\em There exists an interval of values of $a$ (which depends on $p$) inside which the $(p,p,p,\rho)$-regular orbit exists for every $\rho$ in an interval containing $p$. Given $(T_1,T_2)$ the number $\rho$ is unique and its dependence on $(T_1,T_2)$ is as in Theorem \ref{STATEREG}.} \footnote{We have checked the existence conditions of families of $(p,p,p,\rho)$-regular orbit for $p=2,3,4$ and $5$. (The existence conditions for $p=1$ are  given in Theorem \ref{STATEREG}.) For $p\in \{2,3,4,5\}$, the interval of values of $a$ has the form $(a_1(p),a_2(p))$ with $0<a_1(p)<a_2(p)<1$. Notice that for $\rho=p$, the $(p,p,p,\rho)$-regular orbit is the $p$-periodic balanced orbit.}

Along the same lines, the existence of regular orbits other than $(n_A,n_B,n_C,\rho)$-regular ones is unknown. For instance, we do not know if there can be regular orbits which pass $\rho_1\in\Q\setminus\N$ iteration(s) per winding in some atom and $\rho_2\in\Q\setminus\N$ iteration(s) per winding in another atom (regular orbits with non-integer repetitions in two atoms). 
\bigskip
 
\noindent
{\bf Non-regular orbits}

\noindent
In addition to regular orbits, the negative 2-circuit may have orbits in the attractor for which the code cannot be interpreted as given by a rigid rotation on a circle composed of 4 arcs (non regular orbits). For instance for $(a,T_1,T_2)=(0.68,0.7,0.5)$ (a point not filled by the domains of existence of
$(p,\rho)$-regular orbits), there exists a periodic orbit with code
$(10\ 00^2\ 01^3\ 11^2\ 10\ 00\ 01^2\ 11^2)^\infty$ which has period 14. For $(a,T_1,T_2)=(0.1,0.900995,0.9005)$ and for $(a,T_1,T_2)=(0.6,0.58,0.5)$, the  (simplest non regular) periodic code $(10\ 00\ 01\ 11\ 10\ 00\ 01^2\ 11^2)$, which has period 10, is admissible. 

Numerical simulations show however that the typical situation when $T_1\neq T_2$ is a unique orbit in the attractor, although there are cases with several orbits (see e.g. results on balanced orbits, especially Figure \ref{PARADOM}).  
\vskip 1truecm

\centerline{\sc \large Part B. Mathematical analysis}

\section{Analysis of general properties}\label{M-ANAL}
\subsection{The attractor, the global orbits and the admissibility
  condition}\label{M-ATTRA} 
From a mathematical point of view, relation (\ref{MODEL}) is
interpreted as a discrete time dynamical system in $\R^N$ generated by
a piecewise affine contracting map $F$. The object of primary
interested in a dissipative dynamical system is the attractor. The
attractor of $(\R^N,F)$, say $A$, is the largest (forward)
invariant set\footnote{In all the paper, invariant always means {\em
forward invariant}: a set $S$ is said to be invariant if
  $F(S)\subset S$.} for which there exists a bounded neighbourhood
$U\supset A$ so that 
\[
A=\bigcap_{t=0}^{+\infty}F^t(U),
\]
(see \cite{B97} for a discussion on various definitions of attractor). 

Our first statement, which is a consequence of positivity and
normalisation of the interaction weights $K_{ij}$, shows that every
orbit asymptotically approaches the cube $[0,1]^N$. 
\begin{Pro}
$A\subset [0,1]^N$.
\label{LOCAT}
\end{Pro}
{\em Proof:} In order to prove the inclusion, it suffices to show that
$[0,1]^N$ is absorbing, i.e.\ that the image of any ball of radius
$\delta$ around $[0,1]^N$ is included in any smaller ball after a
(sufficiently large) finite number of iterations.

\noindent
The conditions $K_{ij}\geq 0$ and ${\displaystyle\sum_{j\in I(i)}}
K_{ij}=1$ imply the inequalities $a x_i^t\leq x_i^{t+1}\leq a
x_i^t+1-a$ for any $x^t\in\R^N$. (In particular, the inclusion
$F([0,1]^N)\subset [0,1]^N$ follows.) These inequalities imply that
$d(x^{t+1},[0,1]^N)\leq a d(x^t,[0,1]^N)$ where the distance 
$d(\cdot,\cdot)$ is induced by the norm $\|x\|=\max_i|x_i|$ for any $x\in\R^N$. By induction, we obtain
\[
F^t(B_{[0,1]^N}(\delta))\subset B_{[0,1]^N}(a^t\delta)
\]
where $B_{[0,1]^N}(\delta)=\left\{x\in \R^N\ :\
  d(x,[0,1]^N)<\delta\right\}$. \hfill $\Box$
\bigskip

An efficient method to describe the attractor is by using symbolic
dynamics. To that goal we need to consider first the set of points
whose orbit is global, namely $G$. The set $G$ is the set of points
$x\in\R^N$ for which there exists a sequence $\{x^t\}_{t\in\Z}$,
called a global orbit, so that $x^0=x$, $x^{t+1}=F(x^t)$ for all
$t\in\Z$ and $\sup_{t\in\Z}\|x^t\|<+\infty$.

In short terms, points in $G$ are those points with infinite and
bounded past history. Similarly as in \cite{CF97} one proves that a
sequence $\{x^t\}_{t\in\Z}$ is a global orbit iff its components write
\begin{equation}
x_i^t=(1-a)\sum_{k=0}^{+\infty}a^k\sum_{j\in I(i)}K_{ij}\theta_{ij}^{t-k-1},
\label{GLOBORB}
\end{equation}
and we have $\theta_{ij}^t=H(s_{ij}(x_j^t-T_{ij}))$ for all $i=1,N$,
$j\in I(i)$ and $t\in\Z$. In other words every global orbit is
entirely characterised by its code.

Consequently, in order to ensure that the attractor can be described by using
symbolic dynamics, it suffices to show that it coincides with $G$. This
is the scope of the next statement. 
\begin{Pro}
The attractor $A$ and the set $G$ of global orbit components coincide.
\label{AEQUG}
\end{Pro}
By definition the attractor of an arbitrary dynamical system satisfies
$F(A)\subset A$. This statement shows that in our case, the inclusion
is not strict, i.e.\ $F(A)=A$. 

\noindent
{\em Proof:} Assume that $x^0\in G$ and let $\{x^t\}_{t\in\Z}$ be the
corresponding global orbit. The relation (\ref{GLOBORB}) and
positivity and normalisation of the $K_{ij}$'s imply that $x^t\in
[0,1]^N\subset B_{[0,1]^N}(\delta)$ for all $t\leq 0$ where $\delta>0$
is arbitrary, i.e.\ 
\[
x^0\in\bigcap_{t=0}^{+\infty}F^t(B_{[0,1]^N}(\delta)),
\] 
and then $G\subset\bigcap_{t=0}^{+\infty}F^t(B_{[0,1]^N}(\delta))$ for any
$\delta>0$.  

\noindent
In the proof of Proposition \ref{LOCAT}, we have shown that
$F(B_{[0,1]^N}(\delta))\subset B_{[0,1]^N}(\delta)$. Consequently, the
set $\bigcap_{t=0}^{+\infty}F^t(B_{[0,1]^N}(\delta))$ is invariant and
is contained in the bounded ball $B_{[0,1]^N}(\delta)$. By definition
of $A$, we conclude that
$\bigcap_{t=0}^{+\infty}F^t(B_{[0,1]^N}(\delta))\subset A$ for every
$\delta>0$ and in particular that $G\subset A$.

\noindent
In addition, let $U$ be the bounded neighbourhood involved in the
definition of $A$, i.e.\ $A=\bigcap_{t=0}^{+\infty}F^t(U)$. Since $U$
is bounded, there exists $\delta'>0$ such that $U\subset
B_{[0,1]^N}(\delta')$. Thus $A\subset
\bigcap_{t=0}^{+\infty}F^t(B_{[0,1]^N}(\delta'))$ and by using the
previous inclusion, we conclude that
$A=\bigcap_{t=0}^{+\infty}F^t(B_{[0,1]^N}(\delta'))$ for some
$\delta'>0$. 

In order to prove that $A\subset G$, we first show that every point
$x\in A$ has a pre-image in $A$. The previous relation shows that, if
$x\in A$, then $x\in F^{t+1}(B_{[0,1]^N}(\delta'))$ for any $t\geq
1$. That is to say, for every $t\geq 1$, there exists $y_t\in
B_{[0,1]^N}(\delta')$ such that $F(F^t(y_t))=x$. 
Given $t\geq 1$, let $z_t=F^t(y_t)$. For every $t\geq 1$, we have $z_t\in
F^{-1}(x)\cap F^t(B_{[0,1]^N}(\delta'))$.

\noindent
When $a>0$, every point in $\R^N$ has a finite and bounded number of pre-images by $F$ (at most one pre-image for each realisation of symbols $\{\theta_{ij}\}$). When $a=0$, the sets $F^t(B_{[0,1]^N}(\delta'))$ ($t\geq 1$) themselves are finite with bounded cardinality.

\noindent
Therefore in both cases, there exists a pre-image $z\in F^{-1}(x)$ such that
$z_{t_k}=z$ for every $k\geq 0$ where $\{t_k\}$ is a strictly
increasing sequence. In other words, we have
$z\in\bigcap_{k=0}^{+\infty}F^{t_k}(B_{[0,1]^N}(\delta'))$ and by
invariance of the ball $B_{[0,1]^N}(\delta')$, we conclude that the
pre-image $z$ of $x$ belongs to $A$. 

\noindent
Now by induction, for every $x\in A$, one constructs a sequence
$\{x^t\}_{t\leq 0}$ such that $x^0=x$,  $x^{t+1}=F(x^t)$ and $x^t\in
A$ for all $t\leq -1$. In other words, $x$ has an
infinite and bounded past history. So $x\in G$ and then $A\subset
G$. \hfill $\Box$ 
\bigskip

Proposition \ref{AEQUG} implies that an orbit belongs to the attractor iff its components
$\{x_i^t\}$ are given by (\ref{GLOBORB}) and
$\theta_{ij}^t=H(s_{ij}(x_j^t-T_{ij}))$ for all $i,j$ and
$t$. Equivalently, a symbolic sequence $\{\theta_{ij}\}_{t\in\Z}$
codes an orbit in $A$ iff the numbers $\{x_i^t\}$ computed with
(\ref{GLOBORB}) satisfy $\theta_{ij}^t=H(s_{ij}(x_j^t-T_{ij}))$ for
all $i,j$ and $t$. Determining the attractor thus amounts to
determining the set of symbolic sequences which satisfy this
condition. It is not difficult to see that this admissibility
condition is equivalent to the following one.

\noindent
{\em Admissibility condition for a symbolic sequence}
$\{\theta_{ij}^t\}_{t\in\Z}$: The numbers computed by using
(\ref{GLOBORB}) satisfy
\begin{equation}
\begin{array}{rccclcl}
{\displaystyle\sup_{t\in\Z\ :\ \theta_{ij}^t=0}}x_j^t&\lesssim& T_{ij}&\leq& 
{\displaystyle\inf_{t\in\Z\ :\ \theta_{ij}^t=1}}x_j^t&\text{if}&s_{ij}=+1\\
{\displaystyle\sup_{t\in\Z\ :\ \theta_{ij}^t=1}}x_j^t&\leq &T_{ij}&\lesssim& 
{\displaystyle\inf_{t\in\Z\ :\ \theta_{ij}^t=0}}x_j^t&\text{if}&s_{ij}=-1
\end{array}\quad i=1,N,\ j\in I(i)
\label{ADMICOND}
\end{equation}
where $\lesssim$ means $<$ if the corresponding bound is attained (the
supremum is a maximum, the infimum a minimum) and means $\leq$
otherwise. 

\subsection{Normalisation of parameters}\label{M-NORMAL}
In section \ref{R-MODEL} the interaction weights have been normalised
${\displaystyle\sum_{j\in I(i)}} K_{ij}=1$ for every $i$. This
simplifying assumption can be relaxed without modifying the dynamics
(up to dilations).

Indeed, assume that the sequence $\{x_i^t\}$ is an orbit of $F$ (not
necessarily in the attractor) for some given parameters $\{K_{ij}\}$
and $\{T_{ij}\}$ (where the weights need not be normalised). Then for
any vector $\{\alpha_i\}_{i=1,N}$ with $\alpha_i\neq 0$, the sequence
$\{\alpha_ix_i^t\}$ is an orbit of $F$ with parameters
$\{\alpha_iK_{ij}\}$ and $\{\alpha_jT_{ij}\}$. By choosing
$\alpha_i=\frac{1}{\sum_{j\in I(i)} K_{ij}}$ (which is always possible
because $\sum_{j\in I(i)} K_{ij}>0$ for every $i$), this new orbit
becomes an orbit of a mapping $F$ with the weights satisfying
${\displaystyle\sum_{j\in I(i)}} K_{ij}=1$ for every $i$. 
\bigskip

Another assumption in section \ref{R-MODEL} is that the interaction
thresholds $T_{ij}$ all belong to $(0,1)$. This assumption is
justified by the following result on the dynamics of circuits.
\begin{Lem}
Assume that the network is a $N$-circuit. If either $T_i<0$ or $T_i>1$ for some $i\in\Z/N\Z$, then the attractor of $F$
consists of a unique fixed point. 
\end{Lem}
{\em Proof:} Let $i$ be such that $T_i\not\in [0,1]$. We have shown in
the proof of Proposition \ref{LOCAT} that the distance between $x^t$
and $[0,1]^N$ goes to 0 when $t$ increases. A simple reasoning proves that 
the symbol associated with $x_i^t$ does not depend on $t$ provided that $t$ is sufficiently large (i.e.\ there exists $\theta_i\in\{0,1\}$ such that, for any initial condition $x^0$, there exists $t'\in\N$ such that
$H(s_i(x_i^t-T_i))=\theta_i$ for all $t\geq t'$). According to the
relation (\ref{CIRCUIT}), we have for all $t\geq t'$ 
\[
x_{i+1}^{t+1}=ax_{i+1}^t+(1-a)\theta_i,
\]
which implies that, for any initial condition, $x_{i+1}^t$ converges monotonically to $\theta_i$.

Therefore the symbol $H(s_{i+1}(x_{i+1}^t-T_{i+1}))$ remains constant, say equals $\theta_{i+1}$, when $t$ is sufficiently large. As before this implies that $x_{i+2}^t$ converges monotonically to $\theta_{i+1}$ independently on the initial condition. By
repeating the argument, one easily proves that every initial condition
converges to a unique fixed point. \hfill $\Box$ 

\subsection{Symmetries}\label{M-SYMET}
As announced in Section \ref{R-SYMET}, equivalences of dynamics between
distinct networks, and between distinct parameter values within a given network, occur which allow to reduce the number of situations to be analysed. These
symmetries follow essentially from the next statement.
\begin{Lem}
Let $F$ be a mapping given by an interaction graph with $N$ genes, and
by parameters $\{s_{ij}\},\{K_{ij}\}$ and $\{T_{ij}\}$, and let $k$ be
fixed. The sequence $\{x_i^t\}_{t\in\N}$ is an orbit of $F$
with $x_j^t\neq T_{ij}$ for every $i,j,t$ iff the sequence
$\{\overline{x}_i^t\}_{t\in\N}$ defined by
\[
\overline{x}_i^t=\left\{\begin{array}{cl}
x_i^t&\text{if}\ i\neq k\\
1-x_k^t&\text{if}\ i=k
\end{array}\right.
\]
is an orbit of $\overline{F}$ with $\overline{x}_j^t\neq
\overline{T}_{ij}$ for every $i,j,t$ where the parameters of
$\overline{F}$ are defined by 
\[
\overline{s}_{ij}=\left\{\begin{array}{cl}
s_{ij}&\text{if}\ j\in I(i)\ \text{and}\ j\neq k\\
-s_{ik}&\text{if}\ k\in I(i)\ \text{and}\ j=k\ 
\end{array}\right.\ \text{if}\ i\neq k\quad\text{and}\quad
\overline{s}_{kj}=\left\{\begin{array}{cl}
-s_{kj}&\text{if}\ j\in I(k)\ \text{and}\ j\neq k\\
s_{kk}&\text{if}\ k\in I(k)\ \text{and}\ j=k
\end{array}\right.
\]
by $\overline{K}_{ij}=K_{ij}$ for every $i,j$ and by
\[
\overline{T}_{ij}=\left\{\begin{array}{cl}
T_{ij}&\text{if}\ j\in I(i)\ \text{and}\ j\neq k\\
1-T_{ik}&\text{if}\ k\in I(i)\ \text{and}\ j=k
\end{array}\right.
\]
\label{LemFLIP}
\end{Lem}
As indicated in Section \ref{R-SYMET}, an example of a network corresponding to $\overline{F}$ is
given in Figure \ref{FIGFLIP}. 

\noindent
{\em Proof:} We check that $\{\overline{x}_i^t\}_{t\in\N}$ satisfies
the induction induced by $\overline{F}$. The relation
$H(s(x-T))=H(-s((1-x)-(1-T)))$ implies that for any $i\neq k$ such that
$k\in I(i)$, we have
$H(s_{ik}(x_k^t-T_{ik}))=
H(\overline{s}_{ik}(\overline{x}_k^t-\overline{T}_{ik}))$.
Hence for any $i\neq k$, we have 
\[
\overline{x}_i^{t+1}=a\overline{x}_i^t+(1-a)\sum_{j\in I(i)}
\overline{K}_{ij}
H(\overline{s}_{ij}(\overline{x}_j^t-\overline{T}_{ij})). 
\]
Moreover the weight normalisation and the relation
$1-H(s(x-T))=H(-s(x-T))$ which holds for all $x\neq T$ imply  
\[
\overline{x}_k^{t+1}=a\overline{x}_k^t+(1-a)\sum_{j\in I(k)}
\overline{K}_{kj}
H(\overline{s}_{kj}(\overline{x}_j^t-\overline{T}_{kj})). 
\]
The Lemma is proved. \hfill $\Box$

When the network is a circuit, and by simplifying notations as in section \ref{R-SYMET}, the parameters of $\overline{F}$ write
\[
\overline{s}_i=\left\{\begin{array}{cl}
s_i&\text{if}\ i\neq k-1,k\\
-s_i&\text{if}\ i=k-1,k
\end{array}\right.
\quad\text{and}\quad
\overline{T}_i=\left\{\begin{array}{cl}
T_i&\text{if}\ i\neq k\\
1-T_k&\text{if}\ i=k
\end{array}\right.
\]
By applying repeatedly Lemma \ref{LemFLIP} so as to maximise the
number of genes for which both incoming and outgoing arrows are
activations, one shows that, with the exception of orbits on
discontinuities, the dynamics of every circuit with $N$ genes is
equivalent to 
\begin{itemize}
\item{} either the dynamics of the circuit where all interactions signs
  are positive (positive circuit),
\item{} or the dynamics of the circuit where all, excepted $s_N$,
 interactions signs are positive (negative circuit).
\end{itemize}

In order to exhibit the internal symmetry $S$ with respect to the
centre of $[0,1]^N$, starting from the
parameter vectors $(s_1,s_2,\cdots,s_N)$ and
$(T_1,T_2,\cdots,T_N)$, one applies Lemma \ref{LemFLIP} with
$k=1$. This results in a circuit with parameter vectors $(-s_1,s_2,\cdots,-s_N)$ and
$(1-T_1,T_2,\cdots,T_N)$. Then applying Lemma \ref{LemFLIP} with $k=2$ produces a circuit with parameter $(s_1,-s_2,\cdots,-s_N)$ and $(1-T_1,1-T_2,\cdots,T_N)$. By repeating the process until $k=N$ results in a circuit with parameters  $(s_1,s_2,\cdots,s_N)$ and $(1-T_1,1-T_2,\cdots,1-T_N)$ and the desired
symmetry follows (see section \ref{R-SYMET2}).

Finally, the symmetry permutation $R$ is obvious in a circuit with all
signs being positive. In a circuit where all signs but $s_N$ are
positive, one has to combine the permutation with Lemma \ref{LemFLIP}. with $i=1$ in
order to preserve the signs. The result is the $\sigma\circ R$
symmetry where
\[
(\sigma\circ R)_i(x)=\left\{\begin{array}{cl}
x_{i-1}&\text{if}\ i\neq 1\\
1-x_N&\text{if}\ i=1
\end{array}\right.
\]

\section{Negative circuit analysis}
This section presents the dynamical analysis of negative feedback circuits with 1 and 2 nodes respectively. As said in section \ref{R-CONTRAC}, the dynamics of the self-inhibitor essentially follows from results on piecewise affine contracting rotations. On the other hand the analysis and results of the negative 2-circuit are fully original.
  
The analysis of the positive circuit with 2 nodes relies on results on the self-activator and is left to the reader. (The analysis of the self-activator is trivial.)
 
\subsection{The piecewise affine contracting rotation}\label{M-CONTRAC}
The simplest circuit which requires a proper analysis is the circuit
with one node and negative self-interaction. The corresponding map
writes $F(x)=ax+(1-a)H(T-x)$ whose asymptotic dynamics takes place in
the invariant absorbing interval $(aT,aT+1-a]$ (see Figure
\ref{PWCONT}). In order to 
investigate this dynamics in a parameter independent interval, we consider the map $\tilde{F}$ defined on $(0,1]$ by $X\circ F\circ X^{-1}$ where the map $X(x):=\frac{x-aT}{1-a}$ maps $(aT,aT+1-a]$ onto
$(0,1]$. 

The map $\tilde{F}$ can be viewed as a piecewise affine contracting
rotation on the circle.  Adapting the analysis developed in
\cite{C99} for such rotations, its dynamics can be entirely described. 
As announced in section \ref{R-CONTRAC}, the basic result is that all
orbits are asymptotically given by a rotation. 
\begin{Thm}
Independently of the parameters $(a,T)$, for any initial condition $x\in (0,1]$, we have 
\[
\lim_{t\to +\infty}\left({\tilde{F}}^t(x)-\phi(\nu t+\alpha)\right)=0
\]
where $\alpha\in\R$ depends on $x$ and where the rotation number $\nu\in [0,1)$ and the function $\phi:\R\to [0,1]$ only depend on the parameters but not on $x$. The function $\phi$ is 1-periodic and its restriction to $(0,1]$ is left continuous and increasing.
\label{CONJUG}
\end{Thm}
In particular, when the rotation number is rational, every orbit is
asymptotically periodic. When the rotation number is irrational, every orbit is
asymptotically quasi-periodic. 

The dependence of the rotation number $(a,T)\mapsto \nu(a,T)$ on parameters is known explicitely \cite{C99} . It can be expressed as follows
\begin{equation}
\nu(a,T)=\nu\quad\text{iff}\quad \overline{T}(a,\nu)\leq
T\leq\overline{T}(a,\nu-0),
\label{INETHRE}
\end{equation}
where\footnote{$\lfloor x\rfloor$ is the largest integer not larger than $x$.}
\[
\overline{T}(a,\nu)=1-\frac{(1-a)^2}{a}\sum_{k=0}^\infty a^k\lfloor 
\nu (k+1)\rfloor.
\]
Indeed, independently of $a$, this expression shows that the map
$\nu\mapsto\overline{T}(a,\nu)$ is right continuous and strictly
decreasing. Hence, given $(a,T)$ the previous inequalities actually define a
unique number $\nu(a,T)$.

\noindent
Moreover and independently of $a$, the map $T\mapsto\nu(a,T)$ is decreasing and
continuous. The properties $\overline{T}(a,1)=0$ and
$\overline{T}(a,0-0)=1$ imply that its range is $[0,1)$. In addition since
$\overline{T}(a,\nu)<\overline{T}(a,\nu-0)$ when $\nu$ is rational,
it has a Devil's staircase structure (see Figure \ref{ROTANUMB}).

\noindent
Although the map $\phi$ also depends on parameters, it remains unchanged when $T$
moves in the interval $[\overline{T}(a,\nu),\overline{T}(a,\nu-0)]$
(given $a$ and $\nu$ fixed). 

\noindent
The attracting sequences
$\{\phi(\nu(a,T)t+\alpha)\}_{t\in\N}$ (where 
$\alpha$ is arbitrary) are orbits of $\tilde{F}$ in most
cases. However, if the interaction threshold $T$ belongs to
the right boundary of an interval with rational velocity, i.e.\ if 
$T=\overline{T}(a,\nu-0)$ for some $\nu\in\Q$, all sequences 
$\{\phi(\nu(a,T)t+\alpha)\}_{t\in\N}$ are translated of the same
periodic ghost orbit. In that case,
the attractor $A$ is empty \cite{C99}. The set
$\{\overline{T}(a,\nu-0)\ :\ \nu\in \Q\}$ is countable and nowhere
dense as can be deduced from the properties of the function 
$\overline{T}(a,\nu)$ stated above.

\noindent
In all other cases, i.e.\ for the interaction threshold not at 
the right boundary of an interval with rational velocity, the asymptotic
dynamics is semi-conjugated to the rotation $R_{\nu(a,T)}(x)=x+\nu(a,T)\
\text{mod}\ 1$ on the unit circle. Strictly speaking we have
$\tilde{F}\circ\phi=\phi\circ R_{\nu(a,T)}$.

\subsection{The negative 2-circuit}\label{M-NEGAT}
In the negative 2-circuit, the rotation symmetry $\sigma\circ R$
applies to orbits not intersecting discontinuities (see section
\ref{R-SYMET2}). In order to extend this symmetry to all orbits 
(and subsequently to simplify the analysis of their admissibility), we consider the map $F$ defined by
\[
F(x_1,x_2)=(ax_1+(1-a)\theta_2,ax_2+(1-a)\theta_1)
\]
where
\[
(\theta_1,\theta_2)=\left\{\begin{array}{ccccc}
(1,0)&\text{if}&x_1>T_1&\text{and}&x_2\geq T_2\\
(0,0)&\text{if}&x_1\leq T_1&\text{and}&x_2>T_2\\
(0,1)&\text{if}&x_1< T_1&\text{and}&x_2\leq T_2\\
(1,1)&\text{if}&x_1\geq T_1&\text{and}&x_2< T_2\\
\end{array}\right.
\]
This map and the map induced by (\ref{CIRCUIT}) are equal everywhere
except on discontinuities. They have the same symbolic graph (Figure
\ref{NEGCIR}) and the same formal expression of orbits in the
attractor (relation (\ref{GLOBORB}) for the negative 2-circuit), namely
\begin{equation}
x_i^t=(1-a)\sum_{k=0}^\infty a^k\theta_{3-i}^{t-k-1}\ \text{for
  all}\ t\in\Z,\ i=1,2.
\label{GLOBNEG}
\end{equation}
The admissibility condition associated with the present map slightly differs
from (\ref{ADMICOND}). It is entirely compatible with the rotation
symmetry; {\sl the admissibility domain in threshold
space of the image of a symbolic sequence by $\sigma\circ R$ is the
image under the same action of the original admissibility domain.}

\subsubsection{Regular codes and families of regular codes}
As reported in section \ref{R-NEGAT}, the recurrence induced by 
the symbolic graph suggests to
consider the symbolic sequences generated by the rotation $x\mapsto x+\nu\
\text{mod}\ 1$ ($\nu>0$) on the unit circle composed of 4 arcs with length respectively $A,B,C$
and $D$ and corresponding to the 4 atoms of the partition (see Figure \ref{REGULAR}). 
Up to a choice of the origin on the circle, these symbolic sequences (called {\bf regular symbolic sequences} or {\bf regular codes}) are given by 
\begin{equation}
(\theta_1^t,\theta_2^t)=\left\{\begin{array}{ccrcl}
(1,0)&\text{if}&0&\leq\nu t-\lfloor\nu t\rfloor<&A\\
(0,0)&\text{if}&A&\leq\nu t-\lfloor\nu t\rfloor<&A+B\\
(0,1)&\text{if}&A+B&\leq\nu t-\lfloor\nu t\rfloor<&A+B+C\\
(1,1)&\text{if}&A+B+C&\leq\nu t-\lfloor\nu t\rfloor<&A+B+C+D=1
\end{array}\right.\ \text{for all}\ t\in\Z
\label{REGCOD}
\end{equation}
i.e.\
\begin{equation}
\theta_1^t=1+\lfloor\nu t-(A+B+C)\rfloor-\lfloor\nu
t-A\rfloor\quad\text{and}\quad 
\theta_2^t=1+\lfloor\nu t-(A+B)\rfloor-\lfloor\nu t\rfloor\ \text{for all}\ t\in\Z
\label{CONDENSED}
\end{equation}
and are denoted by the 5-uple $(A,B,C,D,\nu)$. 

\noindent
The parameters $A,B,C,D$
and $\nu$ must satisfy the following conditions. They are positive real numbers. The
lengths satisfy the circle normalisation $A+B+C+D=1$. Moreover, in order to generate a symbolic sequence compatible with the symbolic graph of Figure \ref{NEGCIR}, the orbit of the rotation $x\mapsto x+\nu\ \text{mod}\ 1$ must pass each arc. Hence, we must have
\[
\nu\leq I\leq 1-3\nu\ \text{for all}\ I\in\{A,B,C,D\},
\]
and
\begin{equation}
2\nu\leq I+J\leq 1-2\nu\ \text{for all}\ I\neq J\in \{A,B,C,D\}
\label{CONSTR}
\end{equation}
which in particular imply $\nu\leq\frac{1}{4}$. 
\bigskip

By analogy with the self-inhibitor, \footnote{{\bf Regular symbolic sequences in the self-inhibitor:} The results in section \ref{M-CONTRAC} imply that the admissible symbolic sequences of the one-dimensional mapping $\tilde{F}$ can be viewed as being generated by the rotation $x\mapsto x+\nu\
\text{mod}\ 1$ on the unit circle composed of 2 arcs with respective length $E$ (for the arc corresponding to the atom 1) and $F=1-E$ (for the arc associated with the atom 0). 

\noindent
Under this point of view, Theorem \ref{CONJUG} in that section implies that for any pair $(a,T)$ there exists a unique $\nu:=\nu(a,T)$ such that the sequence generated by the rotation on the circle with arc lengths $E=1-\nu$ and $F=\nu$, either is admissible or is the code of a ghost orbit. The arc lengths depend on $\nu$ in the following way
\[
I=n_I\nu-\lfloor n_I\nu\rfloor \ \text{for}\ I=E\ \text{and}\ F
\]
where $n_E=-1$, $\lfloor n_E\nu\rfloor =-1$, $n_F=1$ and $\lfloor n_F\nu\rfloor =0$.} we will focus on the admissibility of the regular symbolic sequences for which the arc lengths depend on $\nu$ in the following way. Given $n_A,n_B,n_C,n_D\in\Z$, we have
\begin{equation}
I=n_I\nu-\lfloor n_I\nu\rfloor\ \text{for}\ I=A,B,C,\ \text{and}\ D
\label{NIDEF}
\end{equation}
Such codes are denoted by $(n_A,n_B,n_C,n_D,\nu)$. Similarly as for the self-inhibitor, the rotation number $\nu$ is allowed to
vary in some interval. The novelty is that this interval depends on the numbers $n_A,n_B,n_C,n_D$ because the rotation number must be chosen so that the orbit passes every arc. 

\noindent
Hence, any 4-uple $(n_A,n_B,n_C,n_D)$ generates a family of $(n_A,n_B,n_C,n_D,\nu)$-codes with the parameter $\nu$ varying in an appropriate interval.
Moreover the normalisation $A+B+C+D=1$ holds for any code in the family iff we have
\begin{equation}
n_A+n_B+n_C+n_D=0\quad\text{and}\quad 
\lfloor n_A\nu\rfloor +\lfloor n_B\nu\rfloor +\lfloor n_C\nu\rfloor
+\lfloor n_D\nu\rfloor =-1,
\label{CFAMIL}
\end{equation}
for any $\nu$ in the corresponding interval. 
In the sequel, we always assume that these conditions hold. We shall regard the first condition as determining $n_D$ when $n_A,n_B$ and $n_C$ have been given.
\bigskip

Before entering the admissibility analysis, we notice that the families $(n_A,n_B,n_C,n_D,\nu)$ contain all regular symbolic sequences with rational rotation number. Indeed, given a sequence $(A,B,C,D,\nu)$ with $\nu$ rational, the 4 lengths can be
modified (without affecting the code) so as to satisfy $I=n_I\nu-\lfloor n_I\nu\rfloor$. In addition for any irrational rotation number, the sequences $(n_A,n_B,n_C,n_D,\nu)$ form a dense subset of regular symbolic sequences with that rotation number. Therefore the sequences $(n_A,n_B,n_C,n_D,\nu)$ are generic regular symbolic sequences. 

\noindent
The decomposition of arc lengths into $I=n_I\nu-\lfloor n_I\nu\rfloor$ however may not be unique. Indeed two distinct 5-uples $(n_A,n_B,n_C,n_D,\nu)$ may generate the same regular symbolic sequence. However if $(n_A,n_B,n_C,n_D)\neq (n'_A,n'_B,n'_C,n'_D)$ there exists $\nu$ (arbitrarily small) such that the symbolic sequences are $(n_A,n_B,n_C,n_D,\nu)$ and $(n'_A,n'_B,n'_C,n'_D,\nu)$ are distinct.\footnote{The same comments hold for families of symbolic sequences in the self-inhibitor.}
\bigskip

For the sake of simplicity, we shall only consider those families of $(n_A,n_B,n_C,n_D,\nu)$ for which the condition $\lfloor n_I\nu\rfloor=0$ holds for 3 elements in $\{A,B,C,D\}$, say $A,B$ and $C$, for all $\nu$ in the corresponding interval. 
{\bf In the rest of the paper we shall assume} that $\lfloor n_A\nu\rfloor =\lfloor n_B\nu\rfloor = \lfloor n_C\nu\rfloor =0$ for all $\nu\in (0,\frac{1}{n_A+n_B+n_C+1}]$ (the corresponding interval of $\nu$ in this case). In phase space, it means that the orbit passes (independently of $\nu$) the same number of iterations per winding in 3 atoms, precisely $n_A$ iterations in 10, $n_B$ in 11 and $n_C$ in 01.  

\subsubsection{Reduction of admissibility for the sequences $(n_A,n_B,n_C,n_D,\nu)$}
In the case of the sequences $(n_A,n_B,n_C,n_D,\nu)$, the admissibility condition turns out to reduce to a condition on a real function. To see this, given $z\in\R$, $n\in\Z$ and $\nu>0$, let\footnote{$\lceil x\rceil$ is the smallest integer not
smaller than $x$.}
\begin{equation}
\varphi(z,n,\nu)=\lceil n\nu\rceil+(1-a)\sum_{k=0}^\infty
a^k\left(\lfloor z-(k+n+1)\nu\rfloor-\lfloor z-(k+1)\nu\rfloor\right)
\label{ORIGNAL}
\end{equation}
which is a 1-periodic function of $z$. 

\noindent
Given any integers
$n_\alpha,n_\beta$ and $n_\gamma$, let
$n_{\alpha\beta}=n_\alpha+n_\beta$ and
$n_{\alpha\beta\gamma}=n_{\alpha\beta}+n_\gamma$. 

\noindent
Together with assumption $\lfloor n_A\nu\rfloor =\lfloor n_B\nu\rfloor =0$ the property $\lfloor A+B\rfloor=0$ implies $\lfloor
n_{AB}\nu\rfloor=\lfloor n_A\nu\rfloor+\lfloor 
n_B\nu\rfloor=0$ and thus $A+B=n_{AB}\nu$. Similarly, one shows that $B+C=n_{BC}\nu$, $A+B+C=n_{ABC}\nu$ and, using also the relation (\ref{CFAMIL}), that $B+C+D=1-n_A\nu$. 
\bigskip

According to the definition of $F$ and to the expression
(\ref{REGCOD}), the symbolic sequence
$(n_A,n_B,n_C,n_D,\nu)$ is admissible iff for all $t\in\Z$ we have
\[
\left\{\begin{array}{ccccrcl}
x_1^t>T_1&\text{and}&x_2^t\geq T_2&\text{if}&0&\leq\nu t-\lfloor\nu
t\rfloor<&A\\
x_1^t\leq T_1&\text{and}&x_2^t>T_2&\text{if}&A&\leq\nu t-\lfloor\nu t\rfloor<&
A+B\\
x_1^t< T_1&\text{and}&x_2^t\leq T_2&\text{if}&A+B&\leq\nu t-\lfloor\nu
t\rfloor<&A+B+C\\
x_1^t\geq T_1&\text{and}&x_2^t<T_2&\text{if}&A+B+C&\leq\nu
t-\lfloor\nu t\rfloor<&1 
\end{array}\right.
\]
where $x_1^t$ and $x_2^t$ are computed by inserting the condensed expression (\ref{CONDENSED}) of a regular code into the
expression (\ref{GLOBNEG}) of global orbits. It results from relations (\ref{CONDENSED}) and (\ref{NIDEF}) that the orbit components are given by 
\[
x_1^t=\varphi(\nu t-\lfloor\nu t\rfloor,n_{AB},\nu)\quad \text{and}\quad
x_2^t=\varphi(\nu t-\lfloor\nu t\rfloor-A,n_{BC},\nu)\quad \text{for all}\ t\in\Z
\]
Therefore the symbolic sequence $(n_A,n_B,n_C,n_D,\nu)$ is admissible iff the inequalities (\ref{ADMINEG}) below holds both with $(\alpha,\beta,\gamma,T)=(A,B,C,T_1)$ and with $(\alpha,\beta,\gamma,T)=(B,C,D,T_2)$. The inequalities are
\begin{equation}
\left\{\begin{array}{ccrcl}
\varphi (\nu t-\lfloor\nu t\rfloor,n_{\alpha\beta},\nu)>T&\text{if}&0&\leq \nu t-\lfloor\nu t\rfloor<&n_\alpha\nu-\lfloor
n_\alpha\nu\rfloor\\
\varphi (\nu t-\lfloor\nu t\rfloor,n_{\alpha\beta},\nu)\leq T&\text{if}&n_\alpha\nu-\lfloor
n_\alpha\nu\rfloor &\leq \nu t-\lfloor\nu t\rfloor<&n_{\alpha\beta}\nu-\lfloor
n_{\alpha\beta}\nu\rfloor\\
\varphi (\nu t-\lfloor\nu t\rfloor,n_{\alpha\beta},\nu)<T&\text{if}&n_{\alpha\beta}\nu-\lfloor
n_{\alpha\beta}\nu\rfloor &\leq \nu t-\lfloor\nu t\rfloor<&n_{\alpha\beta\gamma}\nu-\lfloor
n_{\alpha\beta\gamma}\nu\rfloor\\
\varphi(\nu t-\lfloor\nu t\rfloor,n_{\alpha\beta},\nu)\geq T&\text{if}&n_{\alpha\beta\gamma}\nu-\lfloor
n_{\alpha\beta\gamma}\nu\rfloor &\leq \nu t-\lfloor\nu t\rfloor<&1
\end{array}\right.\ \text{for all}\ t\in\Z
\label{ADMINEG}
\end{equation}
Together with the assumption $\lfloor n_A\nu\rfloor =\lfloor n_B\nu\rfloor = \lfloor n_C\nu\rfloor =0$ the next statement allows to simplify further the admissibility condition.
\begin{Lem}
For $\lfloor n_\alpha\nu\rfloor =\lfloor n_\beta\nu\rfloor =0$ and $\nu>0$, the
inequalities (\ref{ADMINEG}) are equivalent to the following ones
\[
\varphi((n_\alpha+1)\nu-0,n_{\alpha\beta},\nu)\lesssim T\leq
\varphi((n_\alpha-1)\nu,n_{\alpha\beta},\nu)
\]
and 
\[
\varphi(n_{\alpha\beta\gamma}\nu-\lfloor n_{\alpha\beta\gamma}\nu\rfloor 
-0,n_{\alpha\beta},\nu)\lesssim T\leq
\varphi(n_{\alpha\beta\gamma}\nu-\lfloor n_{\alpha\beta\gamma}\nu\rfloor
,n_{\alpha\beta},\nu)
\]
where again $\lesssim$ means $<$ when the limits are attained (which happens iff $\nu$ is rational) and means $\leq$ otherwise.
\label{REDUC1}
\end{Lem}
{\sl Proof:} Together with $\nu>0$ the relations (\ref{CONSTR}) and (\ref{NIDEF}) and the conditions $\lfloor
n_\alpha\nu\rfloor =\lfloor n_\beta\nu\rfloor =0$ imply that
$n_\alpha>0,n_\beta>0$ and $\lfloor n_{\alpha\beta}\nu\rfloor =0$. 
When $n_{\alpha\beta}>0$, the expression of $\varphi(z,n_{\alpha\beta},\nu)$ can be rewritten as the sum of two functions, namely $\varphi:=\varphi_1+\varphi_2$ where (note that $\lceil n_{\alpha\beta}\nu\rceil =1$ in the present case)
\[
\varphi_1(z,n_{\alpha\beta},\nu)=1-(1-a)\sum_{k=0}^{n_{\alpha\beta}-1}a^k\lfloor
z-(k+1)\nu\rfloor\quad\text{and}\quad
\varphi_2(z,n_{\alpha\beta},\nu)=(1-a)(a^{-n_{\alpha\beta}}-1)\sum_{k=n_{\alpha\beta}}^\infty a^k\lfloor
z-(k+1)\nu\rfloor.
\]
The map $z\mapsto\varphi_1(z,n_{\alpha\beta},\nu)$ is a right continuous decreasing
function on $[0,1)$ which is constant on every interval $[j\nu,(j+1)\nu)$ where $0\leq j<n_{\alpha\beta}$ (note that $[j\nu,(j+1)\nu)\subset [0,1)$ for any $0\leq j<n_{\alpha\beta}$ because $\lfloor n_{\alpha\beta}\nu\rfloor =0$). It is also constant on the interval $(n_{\alpha\beta}\nu,1]$.

\noindent
The map $z\mapsto\varphi_2(z,n_{\alpha\beta},\nu)$ is a right continuous increasing
function on $[0,1)$. It is a step function if $\nu$ is rational with
discontinuities at $p\nu-\lfloor p\nu\rfloor$ where $p>n_{\alpha\beta}$. When
$\nu$ is irrational, the discontinuities are dense in
$[0,1)$ and the map is strictly increasing. As a consequence, every left
limit $\varphi(z-0,n_{\alpha\beta},\nu)$ is attained when $\nu$ is rational and not
attained when $\nu$ is irrational.

\noindent
Therefore when $\lfloor n_\alpha\nu\rfloor =\lfloor n_\beta\nu\rfloor
=0$ the first two lines in (\ref{ADMINEG}) are equivalent to the following
inequalities
\[
\max_{n_\alpha< p\leq n_{\alpha\beta}}\varphi(p\nu-0,n_{\alpha\beta},\nu)
\lesssim T<
\min_{0\leq p<n_\alpha}\varphi(p\nu,n_{\alpha\beta},\nu)
\]
and since $n_{\alpha\beta}\nu<n_{\alpha\beta\gamma}\nu-\lfloor n_{\alpha\beta\gamma}\nu\rfloor$ (because $\lfloor n_{\alpha\beta\gamma}\nu\rfloor=\lfloor n_{\gamma}\nu\rfloor$ and $0<n_{\gamma}\nu-\lfloor n_{\gamma}\nu\rfloor$), the last two lines in (\ref{ADMINEG}) are equivalent to the following
inequalities
\[
\varphi(n_{\alpha\beta\gamma}\nu-\lfloor n_{\alpha\beta\gamma}\nu\rfloor 
-0,n_{\alpha\beta},\nu)\lesssim T\leq 
\varphi(n_{\alpha\beta\gamma}\nu-\lfloor n_{\alpha\beta\gamma}\nu\rfloor
,n_{\alpha\beta},\nu).
\]
It remains to show that 
\[
\max_{n_\alpha< p\leq n_{\alpha\beta}}\varphi(p\nu-0,n_{\alpha\beta},\nu)=
\varphi((n_\alpha+1)\nu-0,n_{\alpha\beta},\nu)\quad\text{and}\quad
\min_{0\leq p<n_\alpha}\varphi(p\nu,n_{\alpha\beta},\nu)=
\varphi((n_\alpha-1)\nu,n_{\alpha\beta},\nu).
\]
We prove the second assertion; the first one follows similarly. In the case where $n>0$ and $n\nu\in (0,1]$, the expression (\ref{ORIGNAL}) of the function $\varphi$ becomes
\begin{equation}
\varphi(z,n,\nu)=\left\{\begin{array}{clc}
{\displaystyle a^{\lfloor\frac{z}{\nu}\rfloor}-(a^{-n}-1)\sum_{j=1}^{+\infty}a^{\lfloor\frac{z+j}{\nu}\rfloor}}&\text{if}&0\leq z\leq n\nu\\
1-(a^{-n}-1){\displaystyle\sum_{j=0}^{+\infty}}a^{\lfloor\frac{z+j}{\nu}\rfloor}&\text{if}&n\nu\leq z<1
\end{array}\right.
\label{LATEST}
\end{equation}
By assumptions on $n_\alpha,n_\beta$ and $\nu$, we have $\lfloor p\nu\rfloor=0$ and $\lceil (n_{\alpha\beta}-p)\nu\rceil=1$ for every $0\leq p<n_\alpha$. It results that 
\[
\varphi((p+1)\nu,n_{\alpha\beta},\nu)-\varphi(p\nu,n_{\alpha\beta},\nu)=
(1-a)a^p\left((a^{-n_{\alpha\beta}}-1)\psi(\nu)-1\right)
\]
for every $0\leq p<n_\alpha$, where the function $\psi$ is defined in relation (\ref{DEFPSI}) below. 

\noindent
Together with the condition $\lfloor n_{\alpha\beta}\nu\rfloor =0$,
the constraint (\ref{CONSTR}) implies that $\nu\leq\frac{1}{n_{\alpha\beta}+2}$ and we have 
\[
\varphi((p+1)\nu,n_{\alpha\beta},\nu)-\varphi(p\nu,n_{\alpha\beta},\nu)\leq 
(1-a)a^p\left((a^{-n_{\alpha\beta}}-1)\psi(\frac{1}{n_{\alpha\beta}+2})-1\right)<0
\]
since $\psi(\frac{1}{n_{\alpha\beta}+2})=\frac{a^{n_{\alpha\beta}+2}}{1-a^{n_{\alpha\beta}+2}}$. The second assertion is proved. \hfill $\Box$
\bigskip

By computing the quantities of Lemma \ref{REDUC1} both with $(\alpha,\beta,\gamma)=(A,B,C)$ and with $(\alpha,\beta,\gamma)=(B,C,D)$, one can conclude about the existence of a product of intervals in the threshold plane $(T_1,T_2)$ inside which a given sequence $(n_A,n_B,n_C,n_D,\nu)$ is admissible. Such existence may depend on $a$. For an example of such a computation, see the analysis of balanced orbit existence domains in the next section. 

By analogy with the self-inhibitor it is also interesting to determine those triples $(n_A,n_B,n_C)$ -- and those values of $a$ -- for which every sequence $(n_A,n_B,n_C,n_D,\nu)$ with $\nu\in (\nu_1,\nu_2]\subset (0,\frac{1}{n_{ABC}+1}]$ is admissible upon the choice of $(T_1,T_2)$. Indeed, in addition to admissibility, one may also be interested in describing changes of the rotation number with threshold parameters.

\noindent
To that goal the following function $\psi$ is useful 
\begin{equation}
\psi(\nu)=\sum_{j=1}^\infty a^{\lfloor\frac{j}{\nu}\rfloor},\quad \nu>0
\label{DEFPSI}
\end{equation}
This function is strictly increasing, left continuous and with discontinuities for every rational number ($\psi$ is continuous at every irrational number). Moreover, $\psi(0+0)=0$ and the 
relation $\lfloor x-0\rfloor =\lceil x\rceil -1$ implies
$\psi(\nu+0)={\displaystyle\sum_{j=1}^{+\infty}}a^{\lceil\frac{j}{\nu}\rceil-1}$.

\noindent
The conditions on parameters, which ensure that any symbolic sequence in $\{(n_A,n_B,n_C,n_D,\nu)\ :\ \nu\in (\nu_1,\nu_2]\}$ is admissible provided that $(T_1,T_2)$ is suitably chosen, are given the following statement.
\begin{Pro}
Let the integers $n_A,n_B,n_C$, the degradation rate $a\in (0,1)$ and the rotation number interval $(\nu_1,\nu_2]\subset (0,\frac{1}{n_{ABC}+1}]$ be given.

\noindent
There exists a set ${\cal S}$ in threshold space such that, for every $\nu\in (\nu_1,\nu_2]$, the sequence $(n_A,n_B,n_C,n_D,\nu)$ is admissible, upon the choice of $(T_1,T_2)\in {\cal S}$, iff the following conditions hold simultaneously

\noindent
(C1) $(a^{-n_{AB}}-1)(1+a)\psi(\nu_2)\leq 1$ and $(a^{-n_{BC}}-1)(1+a)\psi(\nu_2)\leq 1$.

\noindent
(C2) $0\leq 1-a^{n_A}-a^{n_C}(1-a^{n_{AB}})+a^{n_C}(1-a^{n_{AB}})(a^{-n_{BC}+1}-1)\psi(\nu_1+0)$.

\noindent
(C3) $a^{-n_{ABC}}(1-a^{n_{BC}})(1-a^{n_{AB}+1})\psi(\nu_2)\leq 1-a^{n_B}$.

\noindent
(C4) $(a^{-n_{AB}}-1)(1-a^{n_{BC}+1})\psi(\nu_2)\leq 1-a^{-n_A+1}+a^{n_{BC}}(a^{-n_{AB}}-1)$.

\noindent
(C5) $0\leq 1-a^{-n_B+1}+(a^{-n_{BC}}-1)(a^{-n_{AB}+1}-1)\psi(\nu_1+0)$.
\label{FINALCOND}
\end{Pro}
{\em Proof:} It follows from the previous section that there exists $(T_1,T_2)$ such that the regular symbolic sequence $(n_A,n_B,n_C,n_D,\nu)$ is admissible iff the following conditions hold both with $(\alpha,\beta,\gamma)=(A,B,C)$ and with $(\alpha,\beta,\gamma)=(B,C,D)$
\begin{eqnarray}
0&\lesssim &\Delta_1(n_\alpha,n_\beta,\nu):=\varphi((n_\alpha-1)\nu,n_{\alpha\beta},\nu)-\varphi((n_\alpha+1)\nu-0,n_{\alpha\beta},\nu)\nonumber\\
0&\lesssim &\Delta_2(n_\alpha,n_\beta,n_\gamma,\nu):=\varphi(n_{\alpha\beta\gamma}\nu-\lfloor n_{\alpha\beta\gamma}\nu\rfloor
,n_{\alpha\beta},\nu)-\varphi((n_\alpha+1)\nu-0,n_{\alpha\beta},\nu)\label{ADMISREDUC}\\
0&\lesssim &\Delta_3(n_\alpha,n_\beta,n_\gamma,\nu):=\varphi((n_\alpha-1)\nu,n_{\alpha\beta},\nu)-\varphi(n_{\alpha\beta\gamma}\nu-\lfloor n_{\alpha\beta\gamma}\nu\rfloor 
-0,n_{\alpha\beta},\nu)\nonumber
\end{eqnarray}
where now the symbol $\lesssim$ means $<$ if $\nu$ is rational and $\leq$ if $\nu$ is irrational. (The expected fourth condition, namely 
\[
0\lesssim \varphi(n_{\alpha\beta\gamma}\nu-\lfloor n_{\alpha\beta\gamma}\nu\rfloor
,n_{\alpha\beta},\nu)-\varphi(n_{\alpha\beta\gamma}\nu-\lfloor n_{\alpha\beta\gamma}\nu\rfloor-0,n_{\alpha\beta},\nu)
\]
always holds because $n_{\alpha\beta\gamma}\nu-\lfloor n_{\alpha\beta\gamma}\nu\rfloor$ belongs to the interval $(n_{\alpha\beta}\nu,1]$ where $\varphi$ is increasing, see proof of Lemma \ref{REDUC1}.)

\noindent
The proof thus consists in analysing the conditions (\ref{ADMISREDUC}) with $(\alpha,\beta,\gamma)=(A,B,C)$ and with $(\alpha,\beta,\gamma)=(B,C,D)$ to obtain the desired condition.

\noindent
{\em Analysis of the condition $0\lesssim\Delta_1(n_\alpha,n_\beta,\nu)$.} Similar calculations to those in the proof of Lemma \ref{REDUC1} show that
\[
\Delta_1(n_\alpha,n_\beta,\nu)=a^{n_\alpha -1}\left(1-a-(a^{-n_{\alpha\beta}}-1)(\psi(\nu)-a^2\psi(\nu+0))\right)
\]
and by left continuity of $\psi$
\[
\Delta_1(n_\alpha,n_\beta,\nu-0)=a^{n_\alpha -1}(1-a)\left(1-(a^{-n_{\alpha\beta}}-1)(1+a)\psi(\nu)\right)
\]
Since $n_{\alpha\beta}>1$, we have $a^{-n_{\alpha\beta}}-1>0$. Moreover $\psi(\nu)\lesssim\psi(\nu+0)$ for all $\nu$ and $\Delta_1(n_\alpha,n_\beta,\nu-0)$ is decreasing. It results that the condition $0\leq\Delta_1(n_\alpha,n_\beta,\nu-0)$ for $\nu=\nu_2$ implies (the same condition with $\nu\in(\nu_1,\nu_2]$ and then) $0\lesssim\Delta_1(n_\alpha,n_\beta,\nu)$ for all $\nu\in(\nu_1,\nu_2]$. 

\noindent
On the other hand the condition $0\lesssim\Delta_1(n_\alpha,n_\beta,\nu)$ for all $\nu\in(\nu_1,\nu_2]$ implies, by left continuity of $\psi$, $0\leq \Delta_1(n_\alpha,n_\beta,\nu-0)$ for all $\nu\in(\nu_1,\nu_2]$, i.e.\ $0\leq \Delta_1(n_\alpha,n_\beta,\nu_2-0)$.

\noindent
The conditions in (C1) follow from $0\leq \Delta_1(n_\alpha,n_\beta,\nu_2-0)$ by choosing $(\alpha,\beta)=(A,B)$ and $(\alpha,\beta)=(B,C)$ respectively. 

\noindent
{\em Analysis of the condition $0\lesssim\Delta_2(n_\alpha,n_\beta,n_\gamma,\nu)$.} One shows that 
\[
\Delta_2(n_\alpha,n_\beta,n_\gamma,\nu)=
1-a^{n_\alpha}+a^{n_\gamma}(1-a^{n_{\alpha\beta}})(a^{-n_{\beta\gamma}+1}
\psi(\nu+0)-\psi(\nu))-\left\{\begin{array}{lcl}
a^{n_\gamma}(1-a^{n_{\alpha\beta}})&\text{if}&\lfloor n_\gamma\nu\rfloor =0\\
0&\text{if}&\lfloor n_\gamma\nu\rfloor =-1
\end{array}\right.
\]
For $(\alpha,\beta,\gamma)=(A,B,C)$, we have $\lfloor n_\gamma\nu\rfloor =\lfloor n_C\nu\rfloor =0$. Thus $0\lesssim\Delta_2(n_\alpha,n_\beta,n_\gamma,\nu)$ for all $\nu\in (\nu_1,\nu_2]$ iff 
\[
0\leq 1-a^{n_A}-a^{n_C}(1-a^{n_{AB}})+a^{n_C}(1-a^{n_{AB}})(a^{-n_{BC}+1}-1)\psi(\nu)\quad \nu\in (\nu_1,\nu_2]
\]
from which (C2) immediately follows using that $a^{-n_{BC}+1}-1>0$. 

\noindent
For $(\alpha,\beta,\gamma)=(B,C,D)$, we have $\lfloor n_\gamma\nu\rfloor =\lfloor -n_{ABC}\nu\rfloor =-1$ since we have assumed that $\nu\leq\frac{1}{n_{ABC}+1}$. Thus $0\lesssim\Delta_2(n_\alpha,n_\beta,n_\gamma,\nu)$ for all $\nu\in (\nu_1,\nu_2]$ iff
\[
0\leq 1-a^{n_B}+a^{-n_{ABC}}(1-a^{n_{BC}})(a^{n_{AB}+1}-1)\psi(\nu)\quad\forall \nu\in (\nu_1,\nu_2]
\]
which is equivalent to the condition (C3).

\noindent
{\em Analysis of the condition $0\lesssim\Delta_3(n_\alpha,n_\beta,n_\gamma,\nu)$.} Again similar calculations to those in the proof of Lemma \ref{REDUC1} show that  
\begin{eqnarray*}
\Delta_3(n_\alpha,n_\beta,n_\gamma,\nu)&=&
a^{n_\alpha-1}\left(1-a^{-n_\alpha+1}+(a^{-n_{\alpha\beta}}-1)(a^{n_{
\beta\gamma}+1}\psi(\nu+0)-\psi(\nu))\right)\\
&+&\left\{\begin{array}{lcl}
a^{n_{\alpha\beta\gamma}-1}(a^{-n_{\alpha\beta}}-1)&\text{if}&\lfloor n_\gamma\nu\rfloor =0\\
0&\text{if}&\lfloor n_\gamma\nu\rfloor =-1\end{array}\right.
\end{eqnarray*}
For $(\alpha,\beta,\gamma)=(A,B,C)$, we have $\lfloor n_\gamma\nu\rfloor =0$. Therefore $0\lesssim\Delta_3(n_\alpha,n_\beta,n_\gamma,\nu)$ for all $\nu\in (\nu_1,\nu_2]$ iff
\[
0\leq 1-a^{-n_A+1}+a^{n_{BC}}(a^{-n_{AB}}-1)-(a^{-n_{AB}}-1)(1-a^{n_{BC}+1})\psi(\nu)\quad\forall \nu\in (\nu_1,\nu_2]
\]
which is equivalent to the condition (C4).

\noindent
Finally for $(\alpha,\beta,\gamma)=(B,C,D)$, we have $\lfloor n_\gamma\nu\rfloor =-1$. In this case $0\lesssim\Delta_3(n_\alpha,n_\beta,n_\gamma,\nu)$ for all $\nu\in (\nu_1,\nu_2]$ iff 
\[
0\leq 1-a^{-n_B+1}+(a^{-n_{BC}}-1)(a^{-n_{AB}+1}-1)\psi(\nu)\quad\forall\nu\in (\nu_1,\nu_2]
\]
from which the condition (C5) follows. \hfill $\Box$
\bigskip

Based on Proposition \ref{FINALCOND}, the following statement collects the results on admissibility of a regular symbolic sequence and states the dependence of its rotation number with parameters.
\begin{Pro}
Let the positive integers $n_A,n_B,n_C$, the degradation rate $a\in (0,1)$ and the rotation number interval $(\nu_1,\nu_2]\subset (0,\frac{1}{n_{ABC}+1}]$ be given. Assume that the conditions in Proposition \ref{FINALCOND} hold.

\noindent
Given $\nu\in (\nu_1,\nu_2]$, let $I_1(\nu)$ and $I_2(\nu)$ be the two intervals given by Lemma \ref{REDUC1} such that the regular symbolic sequence $(n_A,n_B,n_C,n_D,\nu)$ is admissible iff $(T_1,T_2)\in I_1(\nu)\times I_2(\nu)$. \footnote{When $\nu$ is irrational, the interval $I_2(\nu)$ reduces to a point.}

\noindent
The boundaries of $I_1(\nu)$ and of $I_2(\nu)$ are strictly decreasing with $\nu$.  In addition, we have $I_2(\nu)<I_2(\nu')$ whenever $\nu'<\nu$ and the union ${\displaystyle\bigcup_{\nu\in (\nu_1,\nu_2]}}I_2(\nu)$ consists of an interval excepted a countable nowhere dense set. 

\noindent
The intervals $I_1(\nu)$ and $I_1(\nu')$ intersect when $\nu$ and $\nu'$ are sufficiently close. 
\label{ENCORE}
\end{Pro}
In particular, given $n_A,n_B,n_C$, when defined, the rotation number is unique and is a decreasing function of $T_2$.

\noindent
{\em Proof:} The left boundary of the interval $I_i(\nu)$ $(i=1,2)$ is the maximum of the quantities involved in the left inequalities of Lemma \ref{REDUC1}. Precisely, the left boundary of $I_1(\nu)$ (resp.\ $I_2(\nu)$) is obtained for  $(\alpha,\beta,\gamma)=(A,B,C)$ (resp.\ $(\alpha,\beta,\gamma)=(B,C,D)$). By using the expression (\ref{LATEST}), these quantities write
\[
\varphi((n_\alpha+1)\nu-0,n_{\alpha\beta},\nu)=a^{n_\alpha}\left(1-a(a^{-n_{\alpha\beta}}-1)\psi(\nu+0)\right)
\]
and
\[
\varphi(n_{\alpha\beta\gamma}\nu-\lfloor n_{\alpha\beta\gamma}\nu\rfloor-0,n_{\alpha\beta},\nu)=1-a^{n_{\alpha\beta\gamma}}(a^{-n_{\alpha\beta}}-1)(a^{-1}(1+\lfloor n_\gamma\nu\rfloor)+\psi(\nu+0))
\]
Both quantities are right continuous strictly decreasing functions of $\nu$. 

\noindent
Similarly, the right boundary of the interval $I_i(\nu)$ is the minimum of the quantities involved in the right inequalities of Lemma \ref{REDUC1}. By using the expression (\ref{LATEST}), these quantities write
\[
\varphi((n_\alpha-1)\nu,n_{\alpha\beta},\nu)=a^{n_\alpha-1}\left(1-(a^{-n_{\alpha\beta}}-1)\psi(\nu)\right)
\]
and 
\[
\varphi(n_{\alpha\beta\gamma}\nu-\lfloor n_{\alpha\beta\gamma}\nu\rfloor,n_{\alpha\beta},\nu)=1-a^{n_{\alpha\beta\gamma}}(a^{-n_{\alpha\beta}}-1)(1+\lfloor n_\gamma\nu\rfloor+\psi(\nu))
\]
They are left continuous strictly decreasing functions of $\nu$.

\noindent
{\sl Dependence of $I_2(\nu)$ on $\nu$:}
The interval $I_2(\nu)$ is obtained for $(\alpha,\beta,\gamma)=(B,C,D)$. In this case, we have $n_{\alpha\beta\gamma}\nu-\lfloor n_{\alpha\beta\gamma}\nu\rfloor=1-n_A\nu$ and $\lfloor n_\gamma\nu\rfloor=-1$.

\noindent
The proof of Proposition \ref{FINALCOND} shows that the condition (C3) implies the following condition: $0\leq\Delta_2(n_B,n_C,n_D,\nu')$ for all $\nu'\in (\nu_1,\nu_2]$. Using the definition of $\Delta_2$ this condition is equivalent to 
\[
\varphi((n_B+1)\nu'-0,n_{BC},\nu')\leq\varphi(1-n_A\nu',n_{BC},\nu')\quad \forall \nu'\in (\nu_1,\nu_2]
\]
Then the right continuity of $\nu'\mapsto\varphi((n_B+1)\nu'-0,n_{BC},\nu')$, the expressions of $\varphi(1-n_A\nu-0,n_{BC},\nu)$ and of $\varphi(1-n_A\nu',n_{BC},\nu')$ and the fact that $\lfloor n_\gamma\nu\rfloor=-1$, imply that for any $\nu\in (\nu_1,\nu_2)$ we have
\[
\varphi((n_B+1)\nu-0,n_{BC},\nu)=\lim_{\nu'\to\nu^+}\varphi((n_B+1)\nu'-0,n_{BC},\nu')\leq\lim_{\nu'\to\nu^+}\varphi(1-n_A\nu',n_{BC},\nu')=\varphi(1-n_A\nu-0,n_{BC},\nu)
\]
Consequently, the left boundary of $I_2(\nu)$ is given by $\varphi(1-n_A\nu-0,n_{BC},\nu)$ for any $\nu\in (\nu_1,\nu_2)$. Similarly, one proves that the right boundary of $I_2(\nu)$ is given by $\varphi(1-n_A\nu,n_{BC},\nu)$ for any $\nu\in (\nu_1,\nu_2]$. By using the expressions above of $\varphi(1-n_A\nu-0,n_{BC},\nu)$ and $\varphi(1-n_A\nu,n_{BC},\nu)$ we conclude that, when $\nu\in (\nu_1,\nu_2)$ is rational, the interval $I_2(\nu)$ is given by
\[
I_2(\nu)=(\xi(\nu+0),\xi(\nu)]
\]
where $\xi(\nu)=\varphi(1-n_A\nu,n_{BC},\nu)=1-a^{-n_A}(a^{-n_{BC}}-1)\psi(\nu)$.\footnote{Moreover if $\nu_2$ is rational, the interval is given by $I_2(\nu_2)=(\max\{\xi(\nu_2+0),\varphi((n_B+1)\nu_2-0,n_{BC},\nu_2)\},\xi(\nu_2)]$.}
In particular, the strict monotonicity of $\psi$ implies that $I_2(\nu)<I_2(\nu')$ whenever $\nu'<\nu$. Moreover, this expression of $I_2(\nu)$ implies that ${\displaystyle\bigcup_{\nu\in (\nu_1,\nu_2]}}I_2(\nu)$ consists of an interval excepted a countable nowhere dense set. 

\noindent
{\sl Dependence of $I_1(\nu)$ on $\nu$:}
The interval $I_1(\nu)$ is obtained for $(\alpha,\beta,\gamma)=(A,B,C)$. In this case, we have $n_{\alpha\beta\gamma}\nu-\lfloor n_{\alpha\beta\gamma}\nu\rfloor=n_{ABC}\nu$ and $\lfloor n_\gamma\nu\rfloor=0$. 

\noindent
By using the expression above of $\varphi(n_{ABC}\nu-0,n_{AB},\nu)$ and of $\varphi(n_{ABC}\nu,n_{AB},\nu)$ we obtain the following inequality
\[
\lim_{\nu'\to\nu^-}\varphi(n_{ABC}\nu'-0,n_{AB},\nu')<\varphi(n_{ABC}\nu,n_{AB},\nu)
\]
In particular, we have $\varphi(n_{ABC}\nu'-0,n_{AB},\nu')<\varphi(n_{ABC}\nu,n_{AB},\nu)$ for all $\nu'<\nu$ sufficiently close.

\noindent
Moreover in the proof of Proposition \ref{FINALCOND} we showed that the left inequality in condition (C1) is equivalent to $0<\Delta_1(n_A,n_B,\nu-0)$ for all $\nu\in (\nu_1,\nu_2)$. By definition of $\Delta_1$, the latter is equivalent to 
\[
\lim_{\nu'\to\nu^-}\varphi((n_A+1)\nu'-0,n_{AB},\nu')<\varphi((n_A-1)\nu,n_{AB},\nu)
\]
which implies that $\varphi((n_A+1)\nu'-0,n_{AB},\nu')<\varphi((n_A-1)\nu,n_{AB},\nu)$ for all $\nu'<\nu$ sufficiently close.

\noindent
Similarly, one shows that condition (C3) and the condition (C4) of Proposition \ref{FINALCOND} imply respectively the following inequalities
\begin{eqnarray*}
\lim_{\nu'\to\nu^-}\varphi((n_A+1)\nu'-0,n_{AB},\nu')&<&\varphi(n_{ABC}\nu,n_{AB},\nu)\\
\lim_{\nu'\to\nu^-}\varphi(n_{ABC}\nu'-0,n_{AB},\nu')&<&\varphi((n_A-1)\nu,n_{AB},\nu)
\end{eqnarray*}
Therefore the interval $I_1(\nu)$ and $I_1(\nu')$ must intersect when $\nu$ and $\nu'$ are sufficiently close. 
\hfill $\Box$

\subsubsection{Examples}\label{EXAMPL-NEGAT}
In this section, we check the admissibility of the sequences $(n_A,n_B,n_C,n_D,\nu)$ in 3 cases. The first case is a family with $\nu\in (0,\nu_2]$ ($\nu_2$ sufficiently close to 0); the second case is a family with $\nu\in (0,\frac{1}{n_{ABC}+1}]$ (complete family) and the third case concerns the admissibility analysis of a sequence with a specific rotation number, namely a balanced orbit.  
\bigskip

\noindent
{\bf 1) Admissibility for every $\nu\in (0,\nu_2]$ -- Proof of Theorem \ref{STATEREG2}.} Theorem \ref{STATEREG2} equivalently states the admissibility of the sequences $(n_A,n_B,n_C,n_D,\nu)$ for every $\nu\in (0,\nu_2]$ (where $\nu_2>0$) depending on threshold parameters. The existence domain and their dependence on the rotation number are given by Proposition \ref{ENCORE}. So we only have to check the conditions (C1)-(C5) of Proposition \ref{FINALCOND}.

\noindent
For $\nu_1=0$, we have $\psi(\nu_1+0)=\psi(0+0)=0$. Together with the assumption $n_B\geq 1$, the condition (C5) with $\nu_1=0$ imposes $n_B=1$. Assuming $\nu_1=0$ and $n_B=1$, the condition (C2) becomes $a^{n_C}\leq\frac{1-a^{n_A}}{1-a^{n_A+1}}$.

\noindent
We now check the conditions involving $\nu_2$. The conditions (C1) and (C3) requires
\[
\psi(\nu_2)\leq\min\{\frac{1}{(a^{-(n_A+1)}-1)(1+a)},\frac{1}{(a^{-(n_C+1)}-1)(1+a)},\frac{a^{n_{AC}+1}(1-a)}{(1-a^{n_C+1})(1-a^{n_A+2})}\}
\]
which holds provided that $\nu_2$ is sufficiently small because the right hand side is positive and $\psi(0+0)=0$. Finally, the condition (C4) imposes that $1-a^{-n_A+1}+a^{n_C+1}(a^{-(n_A+1)}-1)\geq 0$ (i.e.\ that $\frac{a-a^{n_A}}{1-a^{n_A+1}}\leq a^{n_C}$) and is equivalent to 
\[
\psi(\nu_2)\leq\frac{1-a^{-n_A+1}+a^{n_C+1}(a^{-(n_A+1)}-1)}{(a^{-(n_A+1)}-1)(1-a^{n_C+2})}
\]
which holds provided that $\nu_2$ is sufficiently small. Consequently, there exists a positive number $\nu_2$ such that the condition (C1)-(C5) hold with $\nu_1=0$ and that number $\nu_2$ iff  
\[
n_B=1\quad\text{and}\quad \frac{a-a^{n_A}}{1-a^{n_A+1}}\leq a^{n_C}\leq\frac{1-a^{n_A}}{1-a^{n_A+1}}
\]
Analysing these inequalities one concludes that, for any $n_A,n_C\geq 1$, there exist $0\leq\underline{a}_{n_A,n_C}\leq\overline{a}_{n_A,n_C}<1$ such that they hold iff $\underline{a}_{n_A,n_C}\leq a\leq\overline{a}_{n_A,n_C}$. In particular, we have $\underline{a}_{n_A,1}=\underline{a}_{1,n_C}=0$. Theorem \ref{STATEREG2} is proved. 

It may happen that the conditions (C1)-(C5) do no longer hold when $\nu_2$ becomes sufficiently large (in the allowed domain) depending on $a$. However, for $n_A=1$ (or $n_C=1$), the calculations in the next item show that the conditions (C1)-(C5) hold for all $\nu_2$ in the allowed domain provided that $a$ is sufficiently small.
\bigskip

\noindent
{\bf 2) Admissibility of an arbitrary rotation number. Proofs of Theorem \ref{STATEREG} and Proposition \ref{OPTIMAL}.} 
As in the previous proof, since $\nu_1=0$ the condition (C5) implies that an interval of the form $(0,\frac{1}{n_{ABC}+1}]$ can be admissible only if $n_B=1$. 
For the sake of simplicity, we proceed to the analysis in the case where $n_A=1$. The results in the case where $n_C=1$ can be obtained by applying symmetries.  

\noindent
{\em Proof of Theorem \ref{STATEREG}} Theorem \ref{STATEREG} states the existence and the parameter dependence of orbits with code $(1,1,p,-p-2,\nu)$ for all $\nu\in (0,\frac{1}{p+3}]$. Most of this statement is given by Proposition \ref{ENCORE}. Here, we only have to obtain the conditions on $a$ under which every interval $(0,\frac{1}{p+3}]$ is possible.

\noindent
For $n_A=n_B=1$ and $n_C=p$ we have $n_{ABC}+1=p+3$. For $\nu_1=0$ and $\nu_2=\frac{1}{p+3}$, the conditions (C1)-(C4) reduce to (using again $\psi(0+0)=0$ and that $\psi(\frac{1}{p+3})=\frac{a^{p+3}}{1-a^{p+3}}$)
\[
\frac{a^3+a^2+a-1}{a^4+a^2}\leq a^p\leq\frac{1}{1+a}
\]
The condition $a^3+a^2+a\leq 1$ implies that the left inequality holds for any $p\geq 1$. Under the same condition, we have 
\[
a\leq\frac{1}{1+a}
\]
So the previous right inequality also holds for any $p\geq 1$. Therefore, when $a^3+a^2+a\leq 1$, every regular orbit $(1,1,p,-p-2,\nu)$ with $\nu\in (0,\frac{1}{p+3}]$ is admissible provided that $(T_1,T_2)$ is suitably chosen.\hfill $\Box$

\noindent
{\em Proof of Proposition \ref{OPTIMAL}}
The first part of Proposition \ref{OPTIMAL} is a special case of Theorem \ref{STATEREG2}. The second part relies on the fact that, when $a^3+a^2+a-1>0$, then the condition (C3) with $\nu_2=\frac{1}{p+3}$ does not hold for any $p$ sufficiently large. Thus for such $p$, the regular orbits $(1,1,p,-p-2,\nu)$ cannot exists when $\nu$ is close to $\frac{1}{p+3}$.  \hfill $\Box$

\noindent
{\em Example} $a=0.68$ (second picture of Figure \ref{SQUAR_NEG}). In this case,  for any $p\geq 2$, the condition (C3) holds for $\nu_2=\frac{1}{p+4}$ but does not hold when $\nu_2$ is close to $\frac{1}{p+3}$. 
\bigskip

\noindent
{\bf 3) Existence domains of balanced orbits.} The symbolic sequences associated with balanced orbits are regular with $n_A=n_B=n_C=p$ and $\nu=\frac{1}{4p}$. As said in section \ref {R-NEGAT}, the symmetry $\sigma\circ R$ implies that, when non-empty, the corresponding existence domain in the threshold plane is a square centred at $(\frac{1}{2},\frac{1}{2})$ and symmetric with respect to the diagonal. As argued in the beginning of the proof of Proposition \ref{FINALCOND}, such a square exists iff \footnote{Notice that $\Delta_2(p,p,p,\frac{1}{4p})=\Delta_2(p,p,-3p,\frac{1}{4p})$ and $\Delta_3(p,p,p,\frac{1}{4p})=\Delta_3(p,p,-3p,\frac{1}{4p})$.}
\[
0<\Delta_1(p,p,\frac{1}{4p}),\quad 
0<\Delta_2(p,p,p,\frac{1}{4p})\quad\text{and}\quad
0<\Delta_3(p,p,p,\frac{1}{4p})
\]
Computing these quantities explicitely, we find that the first inequality and the second one hold for any $p\geq 1$ and any $a\in [0,1)$. Moreover, the third inequality is equivalent to 
\[
\frac{a^{p-1}}{1+a^{2p}}>\frac{1}{2}
\]
When $p=1$ this inequality holds for any $a$. A simple analysis shows that, for every $p>1$, there exists $\frac{1}{2}<a_p<1$ such that this inequality holds iff $a_p<a<1$. 

\noindent
When this inequality holds, Lemma \ref{REDUC1} indicates that the balanced orbit exists iff
\[
\max\{\varphi(\frac{p+1}{4p}-0,2p,\frac{1}{4p}),\varphi(\frac{3}{4}-0,2p,\frac{1}{4p})\}<T_i<\min\{\varphi(\frac{p-1}{4p},2p,\frac{1}{4p}),\varphi(\frac{3}{4},2p,\frac{1}{4p})\},\quad i=1,2
\]
Focusing on the right bound, one shows that the map $a\mapsto \varphi(\frac{3}{4},2p,\frac{1}{4p})$ is decreasing and tends to $\frac{1}{2}$ when $a\to 1$. Moreover, we have $\varphi(\frac{3}{4},2p,\frac{1}{4p})<\varphi(\frac{p-1}{4p},2p,\frac{1}{4p})$ for every $a$ when $p=1$. When $p>1$, this inequality holds when $a$ is close to 1. When $a$ is close to $a_p$, the converse inequality holds and in this domain the map $a\mapsto \varphi(\frac{p-1}{4p},2p,\frac{1}{4p})$ is increasing, see Figure \ref{PARADOM}.
 
\section{Stability and structural stability of global orbits}\label{M-STABI}
In complement to section \ref{R-RELEV}, the present section contains statements on robustness of orbits not intersecting discontinuities, namely on their Lyapunov and structural stabilities. Lyapunov stability guarantees that the orbit will be observed in an experiment with suitable initial conditions. Structural stability is also relevant because it asserts that the present phenomenology does not depend on the details of the model -- in particular it does not depend on its piecewise affine characteristic. Every orbit not intersecting discontinuities will also be observed in discrete time models generated by smooth (close to piecewise affine) mappings. 
\bigskip

As an extension of the first claim in section \ref{R-RELEV}, assume that $\{x_i^t\}_{t\in\N}$ is an orbit of the system (\ref{MODEL}) at positive distance from discontinuities, i.e.\ we assume that $\inf_{t,i,j\in I(i)}|x_j^t-T_{ij}|>0$. Proposition 2.1 in \cite{LU05} implies that this orbit must be eventually periodic. 
A simple induction shows that this periodic orbit is asymptotically stable. Moreover, its immediate basin of attraction consist of the balls
\[
\{x\in\R^N\ :\ |x_j-x_j^t|<\min_{i\ :\ j\in I(i)}|x_j^t-T_{ij}|,\ \forall j=1,N\}
\]
\bigskip

Now, if $\inf_{t,i,j\in I(i)}|x_j^t-T_{ij}|>0$, then $\inf_{t,i,j\in I(i)}|x_j^t-T'_{ij}|>0$ for any parameters $T'_{ij}$ sufficiently close to $T_{ij}$. An orbit not intersecting discontinuities remains unaffected by small changes in threshold parameters. 

\noindent
Moreover, the expression (\ref{GLOBORB}) depends continuously on $a$ and on $K_{ij}$. As a consequence if, given $(a,K_{ij},T_{ij})$, an orbit satisfies $\inf_{t,i,j\in I(i)}|x_j^t-T_{ij}|>0$, then there exists an orbit in its neighbourhood, with the same symbolic sequence, for any parameters $(a',K'_{ij},T'_{ij})$ sufficiently close to $(a,K_{ij},T_{ij})$. In short terms, every orbit not intersecting discontinuities can be continued for small changes in parameters.
\bigskip

This continuation properties implies a kind of {\sl modularity} of regulation network dynamics, as suggested in \cite{TR99}. Indeed any collection of independent networks, the {\sl units} (which can be circuits but also more complex networks) can be viewed as composing a unique large network in which the weights $K_{ij}$ of interactions between units are zero. By applying the previous continuation argument, one concludes that, in this large network, any combination of unit orbits at positive distance from discontinuities is an orbit at positive distance from discontinuities. Therefore, it can be continued to an orbit of the large network in which all weights are positive - some of the weights being small ({\sl networks composed of weakly interacting units}).
\bigskip

Moreover for an orbit at positive distance from discontinuities, the mapping $F$ ($F^p$ in the case of a $p$-periodic orbit) satisfies the assumptions of the Implicit Function Theorem \cite{Z96}. As a consequence the orbit under consideration has a (unique) continuation for any sufficiently small smooth perturbation of $F$.

\noindent
On the other hand, changes of a mapping in regions of
phase space where a given orbit never enters (e.g.\ close to
discontinuities in the present case) does not affect the
orbit. By combining this argument with the previous one, it results that if an orbit satisfies $\inf_{t,i,j\in I(i)}|x_j^t-T_{ij}|>0$ then it can be continued to any smooth map
sufficiently close to $F$.  

\noindent
Therefore, any of the periodic orbits in the piecewise affine model can be realised in a mapping which involves a non-linear degradation (instead of the linear one) and whose interactions consist of combinations of continuous sigmoidal interactions (instead of the Heaviside function), for instance the Hill function $x\mapsto\frac{x^m}{T^m+x^m}$ (for an activation) or $x\mapsto\frac{T^m}{T^m+x^m}$ (for an inhibition) provided that $m$ is sufficiently large.

\section{Concluding remarks and prospectives}\label{M-LONG}
In this paper, we have introduced an original model for the dynamics of regulatory networks. Just as for other pre-existing models (especially logical networks and ordinary differential equations), the model consists of a dynamical system with several variables representing gene expression levels. The dynamics is based on an interaction graph whose topology reflects the underlying regulatory network structure. The interactions themselves are combinations of sigmoid functions. 

\noindent
The specificity of the present model is to have discrete time and real (continuous) variables. In particular this formalism allows to connect the dynamics of logical networks with the dynamics of ordinary differential equations. Indeed we have shown that both dynamics can be recovered by adjusting a parameter, denoted by $a$, to its extreme values. 

\noindent
This parameter has been argued to quantify interaction delays which can be attributed to finite chemical reaction rates. In other terms, we have shown that a simple discrete time dynamical system can provide qualitative and quantitative insights on delay effects on the dynamics of regulatory networks. Naturally these systems cannot be a substitute to more elaborated models of delay (integro-)differential equations whose phase space is infinite dimensional.
\bigskip

Our analysis of simplest circuits has shown that, just as for logical networks and for ordinary differential equations, any orbit asymptotically approaches a stable fixed point or an oscillatory orbit (stable periodic orbit or a quasi-periodic orbit). Fixed points exist independently of parameters in positive circuits and do not exist in negative circuits (see \cite{S03} for related recent results in ordinary differential equations). 

\noindent
Moreover delays have an effect on the dynamics of the negative self-inhibitor, of the positive 2-circuit and of the negative 2-circuit. These systems have permanent oscillations (depending on initial conditions for the positive 2-circuit). The regular oscillations can be characterised by a rotation number (together with integer repetitions in atoms for the negative 2-circuit) which measures the average time spent in a given region of phase space (atom). 
\bigskip

In all cases, the dependence of the rotation number on parameters has
three characteristic properties which we recall below. Due to finite
accuracy in experiments, when predicting experimental results these
properties may become more important than the exact knowledge of the
rotation number. Once again, they are expected to appear in genetic
regulatory networks with non-negligible interaction delays (and with
degradation rates not (or slightly) depending on genes).

\noindent
The first characteristic property is continuity. The rotation number depends continuously on parameters. It means that small errors on parameters have few impact on the rotation number. 

\noindent
The second property is that the rotation number is constant over intervals of parameters (plateaus). Not only small parameter changes have few impact on the rotation number, but they are likely to keep it unchanged. 

\noindent
The last characteristic property is monotonicity with threshold parameters. Namely, when a threshold parameter increases the average time spent in a given region of phase space increases (or decreases when the rotation number is a decreasing function of the corresponding threshold parameter). Monotonicity is a consequence of the monotonicity of interaction sigmoid functions.
\bigskip

\noindent
{\bf Regular orbits in the negative $N$-circuit}

\noindent
To conclude this paper, we mention that the analysis of regular orbits presented in section \ref{M-NEGAT} can be extended to a negative circuit with arbitrary number of genes. As a representative of negative circuit, we consider the $N$-circuit with all signs equal to 1, excepted $s_N=-1$ (see section \ref{R-SYMET2}).
\begin{figure}
\begin{center}
\input{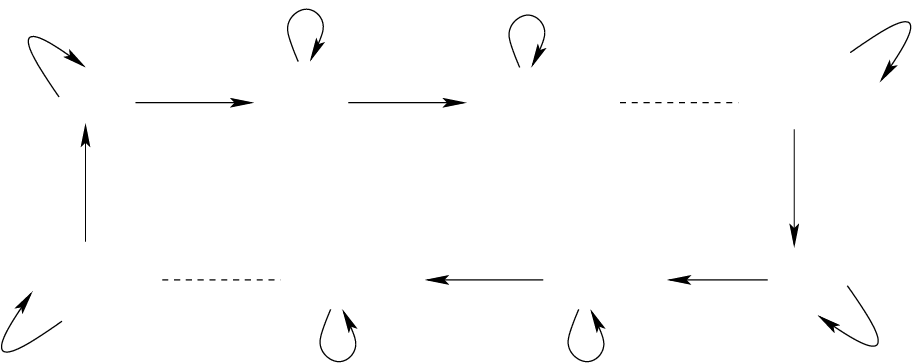}
\end{center}
\caption{Part of the symbolic graph associated with the negative $N$-circuit. The nodes represent the symbol vectors $\{\theta_i\}_{i=1,N}$.}
\label{N-NEGAT}
\end{figure}

\noindent
The corresponding symbolic graph contains the subgraph displayed Figure \ref{N-NEGAT} (which for $N=2$ reduces to the graph of Figure \ref{NEGCIR} (b)). This subgraph (in which any transition from one symbol to the following one only affects one elementary symbol) is the only closed subgraph for which no path can escape  \cite{RMCT03}. (This subgraph is the "attractor" of the asynchronous dynamics in \cite{RMCT03}.) According to the principle of a minimum number of symbol flips during iterations, the orbits corresponding to paths in this subgraph are the most likely to be observed in experiments. 

\noindent
As in the case $N=2$, for any $N\in\N$ one shows that no orbit can stay forever in the same atom and we have a recurrence. As in section \ref{M-NEGAT} the simplest recurrence is a regular orbit (a balanced orbit when the sojourn time does not depend on the atom). 

\noindent
A regular orbit of the negative $N$-circuit is an orbit for which the code is generated by a rotation on the unit circle composed of $2N$ arcs, each arc being associated with a vector symbol (atom) in the subgraph. By assuming that the arc lengths are given as in relation (\ref{NIDEF}) with appropriate normalisation conditions on the integers $n_I$ and $\lfloor n_I\nu\rfloor$, one can proceed to the analysis of the corresponding admissibility condition to conclude as in Theorem \ref{STATEREG} and in Theorem \ref{STATEREG2} about the existence of regular orbits and their changes with parameters.  
\bigskip

\noindent
{\bf Acknowledgements}

\noindent
Using piecewise affine mappings in order to represent genetic
regulatory networks is an idea that grew up at ZiF (University of
Bielefeld) during the programme "The sciences of
complexity". Since fall 2002 when we started investigating this
model, we have benefited from many stimulating discussions (and
emails) and relevant comments from many people, notably A.\ Ciliberto,
C.\ Chandre, C.\ Chaouiya, T.\ Gardner, B.\ Moss\'e, E.\ Remy, D.\
Thieffry and E.\ Ugalde. This work received various financial supports
from, the CNRS, the ECOS-Nord/CONACyT programme M99P01, the Centro
Internacional de Ciencas (UNAM Mexico) and the Grupo de
F\'{\i}sica-Matem\'atica (University of Lisbon).

\end{document}